\documentclass[preprint]{elsarticle} 

\makeatletter
	\def\ps@pprintTitle{%
 	\let\@oddhead\@empty
	\let\@evenhead\@empty
	\def\@oddfoot{\centerline{\thepage}}%
	\let\@evenfoot\@oddfoot}
\makeatother

\usepackage{etoolbox}
\patchcmd{\MaketitleBox}{\footnotesize\itshape\elsaddress\par\vskip36pt}{\footnotesize\itshape\elsaddress\par\parbox[b][36pt]{\linewidth}{\vfill\hfill\textnormal{\today}\hfill\null\vfill}}{}{}%
\patchcmd{\pprintMaketitle}{\footnotesize\itshape\elsaddress\par\vskip36pt}{\footnotesize\itshape\elsaddress\par\parbox[b][36pt]{\linewidth}{\vfill\hfill\textnormal{\today}\hfill\null\vfill}}{}{}%
\usepackage{booktabs}
\usepackage{threeparttable}
\usepackage{algorithm}
\usepackage{algpseudocode}
\usepackage{tikz}
\usetikzlibrary{arrows.meta, positioning, shapes.geometric, shadows}
\usepackage{subcaption}
\usepackage{graphicx}
\usepackage{makecell}
\usepackage{caption}
\usepackage{multirow}
\usepackage{placeins}

\usetikzlibrary{shapes, arrows, positioning, fit, backgrounds}

\usepackage[
            colorlinks=true,%
            breaklinks=true,%
            linkcolor=blue,
            urlcolor=blue,%
            citecolor=blue,%
            pdftitle={Title of paper}, 
            pdfkeywords={Keywords},	
            pdfauthor={Authors},		
            bookmarksopen=false,
            pdfpagemode=UseNone]{hyperref}
            
\usepackage[margin = 1.2in]{geometry}
\usepackage[T1]{fontenc}
\usepackage[english]{babel}
\usepackage[latin1]{inputenc} 
\usepackage{mathtools}
\usepackage{amsmath}
\usepackage{amsfonts}
\usepackage{amsthm}
\usepackage{amssymb}
\usepackage{array}
\usepackage{algorithm} 
\usepackage{algorithmicx}
\usepackage{algpseudocode}
\usepackage{subcaption}
\usepackage{siunitx}

\usepackage{color}
\usepackage{url}
\usepackage{cleveref}
\usepackage{graphicx}
\usepackage{breakcites}
\usepackage{soul} 
\usepackage{enumitem}
\usepackage[title,titletoc,toc]{appendix}

\usepackage[textsize=tiny]{todonotes}

\newtheorem{problem}{Problem}

{\noindent {\textbf{Proof}:} }%
{\hfill $\Box$ \\[1ex] }

\renewcommand{\algorithmicrequire}{\textbf{Input:}}
\renewcommand{\algorithmicensure}{\textbf{Output:}}

\newcommand{\bit}{\begin{itemize}}
\newcommand{\eit}{\end{itemize}}
\newcommand{\ben}{\begin{enumerate}}
\newcommand{\een}{\end{enumerate}}

\DeclareMathOperator*{\argmin}{arg\,min}%
\graphicspath{{../}{./}}

\begin{document}
	
	\begin{frontmatter}
		
		\title{Risk-averse design optimization with CVaR constraints via multifidelity tail-region correction}
		
		
		\author[aff1]{Gyeolhee Lee}
        \ead{rufgml98@hanyang.ac.kr}
        
		\author[aff2]{Dongjin Lee\corref{cor1}}
		\ead{dlee46@hanyang.ac.kr}
		
		\cortext[cor1]{Corresponding author}
		
		\address[aff1]{Department of Automotive Engineering (Automotive-Computer Convergence), Hanyang University, Seoul, South Korea}

        \address[aff2]{Department of Automotive Engineering, Hanyang University, Seoul, South Korea}
        
		\begin{abstract}
           Risk-averse design optimization with conditional value-at-risk (CVaR) constraints requires accurate tail estimates, but repeated high-fidelity (HF) evaluations are costly. We propose a multifidelity (MF) tail-region correction method that uses dimensionally decomposed generalized polynomial chaos expansion (DD-GPCE) as a global surrogate and directs a limited HF budget to responses that most strongly influence CVaR. Covariance of the estimated DD-GPCE coefficients quantifies finite-sample prediction uncertainty and defines a confidence-interval-based tail region. Within this region, a two-stage strategy explores locations with high prediction uncertainty and then exploits locations with large correction-induced shifts. A Tikhonov-regularized residual expansion using the same polynomial basis yields a unified tail-corrected surrogate. For a ten-bar truss, the method produced a design close to the crude-MCS reference with approximately 72--274 times fewer HF evaluations than the considered MF importance sampling configurations. For a suction valve, it reduced mass by 1.05\% and bending and bulging CVaR by 17.68\% and 21.53\%, respectively, while satisfying the valve-opening requirement. The results demonstrate accurate, HF-sample-efficient CVaR-constrained design through localized MF correction.
		\end{abstract}	
		
		
		
\end{frontmatter}

\section{Introduction} \label{sec:intro}


Complex mechanical systems face uncertainty from manufacturing tolerances, material properties, and operating conditions. Their design criteria must therefore address both typical behavior and low-probability, high-impact responses. Conventional safety margins do not directly control failure probability, loss quantiles, or failure severity and can yield either overly conservative or insufficiently protected designs~\cite{Chaudhuri2022,MollerHansson2008}. Design optimization under uncertainty is consequently central to safety and reliability in aerospace, civil, and energy systems~\cite{Chaudhuri2022,royset2025risk,yao2011review,RockafellarRoyset2015,stover2023reliability}.

Robust design optimization (RDO) controls response moments, whereas reliability-based design optimization (RBDO) constrains failure probability~\cite{yao2011review,Aoues2010,MoustaphaSudret2019}. RDO does not directly represent rare tail events, and RBDO does not quantify loss severity after failure~\cite{Chaudhuri2022,RockafellarRoyset2010}. Moreover, sampling rare events with high-fidelity (HF) models can be prohibitively expensive. Safety-critical applications therefore require risk measures that account for both the likelihood and severity of extreme responses~\cite{RockafellarRoyset2015}.

Value-at-Risk (VaR) identifies a prescribed loss quantile, whereas conditional value-at-risk (CVaR) measures the mean loss beyond that threshold~\cite{RockafellarUryasev2000,RockafellarUryasev2002}. CVaR is coherent and integrates readily with optimization, making it suitable as an engineering design objective or constraint~\cite{Artzner1999,KouriSurowiec2016,YangGunzburger2017,Royset2017}. Its reliable estimation, however, requires extensive tail sampling and becomes costly with HF models~\cite{LeeKramer2023}.

Surrogate models such as polynomial chaos expansions (PCEs), Gaussian processes, reduced-order models, and neural networks reduce the cost of repeated HF simulations~\cite{Jakeman2022,LeeKramer2023SMO,chen2012stochastic,Heinkenschloss2018,Heinkenschloss2020,barrera2026statistical}. Generalized PCE (GPCE) constructs multivariate orthonormal polynomials with respect to the joint input distribution, thereby accommodating statistically dependent random inputs without assuming independence~\cite{rahman2018polynomial,lee2020practical}. Dimensionally decomposed GPCE (DD-GPCE) preserves this measure-consistent treatment of dependence while restricting the basis to low-order input interactions, making surrogate construction more tractable for high-dimensional problems~\cite{lee2023high}. However, even a globally accurate DD-GPCE surrogate may retain localized tail errors that bias CVaR estimates~\cite{LeeKramer2023,Heinkenschloss2020,Jakeman2022}.

Adaptive sampling can reduce repeated CVaR-estimation cost. Beiser et al.~\cite{beiser2023adaptive}, for example, increased the sample size used for stochastic-gradient estimates during optimization. That approach adapts sample count rather than selecting realizations from local surrogate uncertainty and does not correct localized errors in surrogate-based CVaR constraints. Such errors can misclassify candidate designs, motivating HF allocation according to each sample's influence on the CVaR constraint.

Multifidelity (MF) modeling provides a natural mechanism for this purpose by combining a computationally inexpensive low-fidelity (LF) model with a limited number of high-fidelity (HF) evaluations~\cite{peherstorfer2018survey}. PCE-based models are particularly compatible with additive multifidelity correction because both the baseline response and its discrepancy can be represented using orthogonal polynomial bases derived from the input probability measure~\cite{ng2012multifidelity}. When the same basis is retained for the baseline and residual models, the corrected response can be written as a unified polynomial expansion and updated through a regularized coefficient correction~\cite{bryson2017all}. This structure allows the baseline model to represent the global response trend while a small number of HF evaluations correct localized discrepancies without requiring a separate surrogate-evaluation procedure. DD-GPCE therefore provides a suitable LF baseline for localized multifidelity correction in high-dimensional stochastic design problems.

In this study, we propose a CVaR-constrained design optimization method that uses localized tail correction for efficient surrogate-based CVaR estimation. We first quantify the epistemic prediction uncertainty introduced by estimating the DD-GPCE coefficients from finite training data. We then use this uncertainty to define a confidence-interval based $\epsilon$-risk region that conservatively retains samples that may affect the CVaR estimate. Within this region, we allocate a limited HF evaluation budget through a two-stage adaptive strategy: predictive variance guides initial exploration, and a prediction-shift criterion guides subsequent exploitation. Using the selected HF samples, we construct a Tikhonov-regularized residual model in the DD-GPCE basis, enabling stable correction even when the number of HF samples is smaller than the number of basis terms. Finally, we combine the LF and residual models into a unified polynomial representation and embed it in the optimization loop. 

The main contributions of this study are as follows:
\begin{itemize}
    \item We quantify finite-sample DD-GPCE prediction uncertainty from the covariance of the estimated expansion coefficients and use it to define a confidence-interval-based $\epsilon$-risk region that retains samples that affect CVaR.

    \item We combine two-stage adaptive HF sampling, guided by predictive variance and weighted prediction shift, with a Tikhonov-regularized residual expansion in the same DD-GPCE basis to produce a unified tail-corrected surrogate.

    \item We embed the tail-corrected surrogate in CVaR-constrained design optimization and assess its CVaR estimation accuracy, design feasibility, and HF evaluation cost.
\end{itemize}

The remainder of this paper is organized as follows. Section~\ref{sec:theoretical_background} presents the CVaR-constrained design formulation and reviews DD-GPCE and additive multifidelity modeling, building on our previous studies of bi-fidelity and multifidelity CVaR estimation using DD-GPCE-based surrogates~\cite{LeeKramer2023SMO,LeeKramer2023}. Section~\ref{sec:mf_for_RADO} introduces the confidence interval-based $\epsilon$-risk region, the two-stage adaptive sampling strategy, the regularized residual correction, and the complete optimization framework. Section~\ref{sec:examples} evaluates the proposed method using the Griewank function, a ten-bar truss structure, and an industrial suction valve system. Finally, Section~\ref{sec:Conclusions} summarizes the principal findings and discusses the limitations and directions for future work.

\section{Theoretical background} \label{sec:theoretical_background}

This section defines the random variables and CVaR-constrained design problem, then summarizes the DD-GPCE surrogate and additive multifidelity update used in the proposed method.

\subsection{Preliminaries} \label{subsec:preliminaries}

Let $\mathbb{N}$, $\mathbb{N}_0$, $\mathbb{R}$, and $\mathbb{R}_0^+$ denote the sets of positive integers, nonnegative integers, real numbers, and nonnegative real numbers, respectively. For any given $N \in \mathbb{N}$, let $\mathbb{R}^N$ denote the $N$-dimensional real vector space and let $\mathbb{A}^N \subseteq \mathbb{R}^N$ denote a bounded subdomain.

\subsection{Design variables and input--output random variables} \label{subsec:variables}

Consider an $N$-dimensional input random vector $\mathbf{X}=(X_1,\ldots,X_N)^{\top}$, whose components $X_1,\ldots,X_N$ represent the uncertain model inputs. Among these input random variables, let $X_{i_1},\ldots,X_{i_M}$ ($1\le i_1 < \ldots<i_M \le N$) be the variables whose mean values are treated as controllable design parameters. These mean parameters are collected in the design vector $\mathbf{d} = (d_1, \dots, d_M)^\top \in \mathcal{D} \subseteq \mathbb{R}^M$ where $d_k$ corresponds to the mean of $X_{i_k}$. The admissible design space is the closed hyperrectangle  $\mathcal{D}=\prod_{k=1}^M[d_{k,L},d_{k,U}]$, where $d_{k,L}$ and $d_{k,U}$ are the lower and upper bounds of $d_k$, respectively. 

For a fixed design $\mathbf{d}\in\mathcal{D}$, consider the probability space $(\Omega_{\mathbf{d}},\mathcal{F}_{\mathbf{d}},\mathrm{P}_{\mathbf{d}})$, where $\Omega_{\mathbf{d}}$ is the sample space, $\mathcal{F}_{\mathbf{d}}$ is the $\sigma$-algebra on $\Omega_{\mathbf{d}}$, and $\mathrm{P}_{\mathbf{d}} : \mathcal{F}_{\mathbf{d}} \to [0,1]$ is the probability measure induced by $\mathbf{d}$. On this space, $\mathbf{X}$ is the measurable mapping $\mathbf{X}:\Omega_{\mathbf{d}}\rightarrow \mathbb{A}^N,$ where $\mathbb{A}^N\subseteq \mathbb{R}^N$ is the set of possible input realizations $\mathbf{x}=(x_1,\ldots,x_N)^{\top}$. The selected input means satisfy $d_k=\mathbb{E}_{\mathbf{d}}[X_{i_k}]$, $k=1,\ldots,M$, where $\mathbb{E}_{\mathbf{d}}[\cdot]$ denotes expectation with respect to $\mathrm{P}_{\mathbf{d}}$. Thus, $\mathbf{d}$ specifies the controllable nominal input values, whereas the variability about these values represents the aleatoric uncertainty inherent in the system.

The joint cumulative distribution function (CDF) of $\mathbf{X}$ is $F_{\mathbf{X}}(\mathbf{x};\mathbf{d}) := \mathrm{P}_{\mathbf{d}}\left(\bigcap_{i=1}^{N} \{X_i \le x_i\}\right)$. Assuming that the distribution of $\mathbf{X}$ is absolutely continuous, its joint probability density function (PDF) is $f_{\mathbf{X}}(\mathbf{x};\mathbf{d}) := \partial^N F_{\mathbf{X}}(\mathbf{x};\mathbf{d})/\partial x_1 \cdots \partial x_N.$ The probability space associated with $\mathbf{X}$ can then be represented as $(\mathbb{A}^N, \mathcal{B}^N, f_{\mathbf{X}}(\mathbf{x};\mathbf{d})\mathrm{d}\mathbf{x}),$ where  $\mathcal{B}^N := \mathcal{B}(\mathbb{A}^N)$ is the Borel $\sigma$-algebra on $\mathbb{A}^N$.

Given a measurable model function $y:\mathbb{A}^N \rightarrow \mathbb{R}$, the corresponding output random variable is $y(\mathbf{X}):\Omega_{\mathbf{d}}\rightarrow\mathbb{R}$. Although $y(\cdot)$ does not explicitly depend on $\mathbf{d}$, the distribution of $y(\mathbf{X})$ varies with $\mathbf{d}$ through the probability law of $\mathbf{X}$. This response may represent, for example, stress, displacement, natural frequency, buckling-related behavior, or a crashworthiness measure. Aleatoric uncertainty propagating through the model may produce extreme responses associated with system failure. 

For a system with $J$ outputs, each scalar output random variable is defined separately as $y_j(\mathbf{X}):\Omega_{\mathbf{d}}\rightarrow\mathbb{R}$, $j=1,\ldots,J$, where $y_j:\mathbb{A}^N\rightarrow\mathbb{R}$ is the $j$th measurable model function. 

\subsection{Risk-averse design optimization} \label{subsec:RADO}

Risk-averse design optimization selects $\mathbf{d}$ to control the tail behavior of the response distribution. Common risk measures for this purpose include Value-at-Risk (VaR) and Conditional Value-at-Risk (CVaR)~\cite{RockafellarRoyset2015}.
At each fixed design vector $\mathbf{d}$, the distribution of $y(\mathbf{X})$ is governed by $\mathrm{P}_{\mathbf{d}}$. We adopt the upper-tail convention, under which larger values of $y(\mathbf{X})$ represent less desirable outcomes. For a risk level $\beta \in (0,1)$, the VaR of $y(\mathbf{X})$ is its $\beta$-quantile
\begin{equation}
    \mathrm{VaR}_{\beta;\mathbf{d}}(y(\mathbf{X})) := \inf \{ t \in \mathbb{R} : \mathrm{P}_{\mathbf{d}}(y(\mathbf{X}) \le t) \ge \beta \}. \label{eq:var}
\end{equation}
If $F_{Y;\mathbf{d}}$ denotes the CDF of $y(\mathbf{X})$, then $\mathrm{VaR}_{\beta;\mathbf{d}}(y(\mathbf{X})) = F_{Y;\mathbf{d}}^{-1}(\beta)$, where $F_{Y;\mathbf{d}}^{-1}$ denotes the generalized inverse of $F_{Y;\mathbf{d}}$. Thus, $\mathrm{VaR}_{0.99;\mathbf{d}}$ marks the threshold at which the upper $1\%$ tail of the response distribution begins.

Although VaR has a direct quantile interpretation, it is not subadditive in general and provides no information about the severity of responses beyond the quantile threshold~\cite{Artzner1999,RockafellarUryasev2000}. CVaR addresses these limitations by accounting for the magnitude of the responses in the upper tail. For a distribution that is continuous at the corresponding quantile, CVaR is defined as
\begin{equation}
    \mathrm{CVaR}_{\beta;\mathbf{d}}(y(\mathbf{X})) = \mathbb{E}_{\mathbf{d}}[y(\mathbf{X}) \mid y(\mathbf{X}) \ge \mathrm{VaR}_{\beta;\mathbf{d}}(y(\mathbf{X}))]. \label{eq:cvar_conditional}
\end{equation}
This definition follows the standard upper-tail interpretation of CVaR for continuous response distributions~\cite{RockafellarUryasev2000,RockafellarUryasev2002}.
Equivalently, 
\begin{equation}
    \mathrm{CVaR}_{\beta;\mathbf{d}}(y(\mathbf{X})) = \mathrm{VaR}_{\beta;\mathbf{d}}(y(\mathbf{X})) + \frac{1}{1-\beta}\mathbb{E}_{\mathbf{d}}\left[\left(y(\mathbf{X}) - \mathrm{VaR}_{\beta;\mathbf{d}}(y(\mathbf{X}))\right)^{+}\right], \label{eq:cvar_plus}
\end{equation}
where $(\cdot)^{+} = \max(\cdot,0)$. This representation shows that CVaR accounts for both the VaR threshold and the expected excess beyond that threshold. 

\begin{figure}
    \centering
    \includegraphics[width=0.8\linewidth]{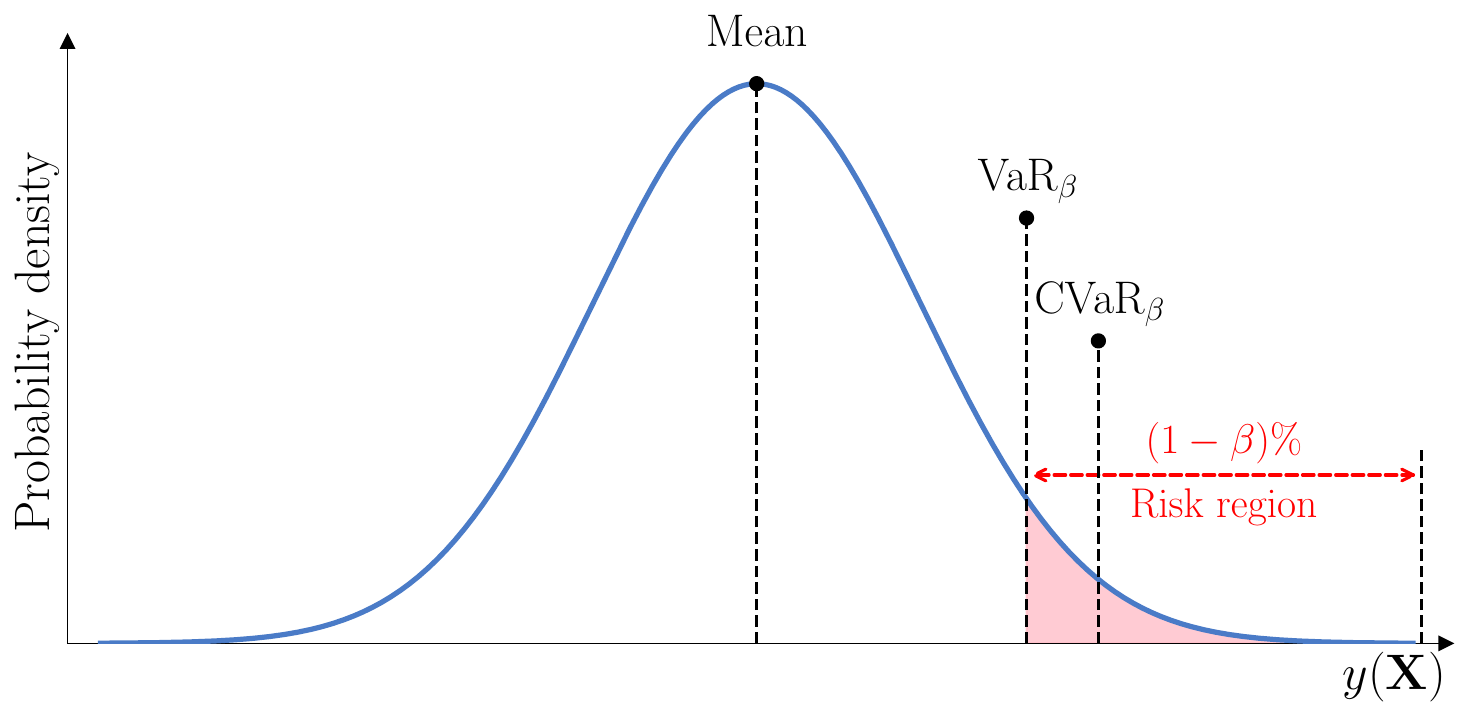}
    \caption{$\mathrm{VaR}_\beta$ and $\mathrm{CVaR}_\beta$ for the upper $100(1-\beta)\%$ risk region.}
    \label{fig:var_cvar_plot}
\end{figure}

Figure~\ref{fig:var_cvar_plot} illustrates the upper-tail region, which contains the worst $100(1-\beta)\%$ of the response distribution. We define the corresponding risk region as
\begin{equation}
    \mathcal{G}_{\beta;\mathbf{d}} := \left\{\mathbf{x} \in \mathbb{A}^N :
    y(\mathbf{x}) \ge \mathrm{VaR}_{\beta;\mathbf{d}}\bigl(y(\mathbf{X})\bigr)
    \right\}. \label{eq:risk_region}
\end{equation}
This region contains the input realizations that contribute to the CVaR evaluation. Under the continuity assumption, CVaR can also be expressed as 
\begin{equation}
    \mathrm{CVaR}_{\beta;\mathbf{d}}(y(\mathbf{X})) = \frac{1}{1-\beta}\mathbb{E}_{\mathbf{d}}\left[y(\mathbf{X})\mathbb{I}_{\mathcal{G}_{\beta;\mathbf{d}}}(\mathbf{X})\right], \label{eq:cvar_indicator}
\end{equation}
where $\mathbb{I}_{\mathcal{G}_{\beta;\mathbf{d}}}(\mathbf{X})$ is the indicator function for the risk region. It captures the magnitude of the tail risk and satisfies the properties of a coherent risk measure, namely monotonicity, subadditivity, positive homogeneity, and translation invariance~\cite{Artzner1999,RockafellarUryasev2002}. When the underlying response functions and feasible set satisfy convexity assumptions, this property can help preserve the convex structure of the resulting optimization problem and facilitate the use of standard convex optimization algorithms~\cite{RockafellarUryasev2000,Chaudhuri2022}.

VaR and CVaR provide complementary descriptions of upper-tail risk. VaR identifies the threshold at which the tail begins, whereas CVaR quantifies the average severity of responses within that tail. A CVaR constraint therefore controls the magnitude of extreme responses rather than only their quantile threshold, making it well suited to engineering systems in which tail-response severity is a primary design concern.

Using CVaR, we formulate the risk-averse design problem as~\cite{RockafellarUryasev2000,KouriSurowiec2016}
\begin{equation}
\begin{aligned}
    \min_{\mathbf{d}\in\mathcal{D}} \quad
    & c(\mathbf{d}) &\\
    \text{subject to} \quad
    & \mathrm{CVaR}_{\beta;\mathbf{d}}
      \!\left(y_j(\mathbf{X})\right)
      \le \tau_j,
      &\qquad j=1,\ldots,J ,\\
      &d_{k,L}\leq d_k \leq d_{k,U}, &\qquad k=1,\ldots,M.
\end{aligned}
\label{eq:multioutput_risk_averse_design}
\end{equation}
Here, \(c(\mathbf{d})\) denotes a design objective, such as mass,
volume, or cost. The function $y_j(\mathbf{X})$ denotes
the $j$th safety-related response, and $\tau_j\in\mathbb{R}$ is its acceptable
risk limit. The responses may represent stress, strain, displacement,
buckling-related behavior, or crashworthiness measures. Thus, the
formulation imposes an individual CVaR constraint on each scalar
response.

Constraining CVaR directly limits the average magnitude of responses within the specified upper tail. Its accurate evaluation, however, remains computationally challenging. The input random vector may be high-dimensional and correlated, the model response $y(\mathbf{X})$ may be strongly nonlinear or nonsmooth, and each model evaluation may require an expensive numerical simulation. Moreover, samples become increasingly sparse in the relevant tail as $\beta$ approaches one, increasing the variance and computational cost of CVaR estimation~\cite{Heinkenschloss2018}. Accurate and efficient risk estimation therefore remains a central challenge in practical risk-averse design optimization~\cite{Heinkenschloss2018,Jakeman2022}. These computational demands motivate the multifidelity surrogate strategy introduced in the subsequent sections.

\begin{problem}
Determine an optimal design that satisfies the prescribed CVaR constraints within a limited HF budget. The central challenge is to reliably estimate tail risk at each candidate design despite sparse tail-region samples and localized surrogate bias.
\end{problem}

\subsection{Dimensionally decomposed generalized polynomial chaos expansion} \label{subsec:ddgpce}

PCE represents a square-integrable response $y(\mathbf{X})$ in an orthogonal polynomial basis~\cite{XiuKarniadakis2002}. Unlike tensor-product PCE for independent inputs, the GPCE used here constructs multivariate polynomials orthonormal with respect to the dependent joint measure $f_{\mathbf{X}}(\mathbf{x};\mathbf{d})\mathrm{d}\mathbf{x}$~\cite{rahman2018polynomial,lee2020practical}. Supplementary Section~S1 provides the full GPCE formulation.

For a total polynomial degree $m$, the full GPCE contains $\binom{N+m}{m}$ basis functions. Its size therefore grows combinatorially with the input dimension $N$ and polynomial degree $m$. In many engineering applications, however, individual input effects and
low-order interactions among a small number of variables dominate the response, whereas higher-order interactions make comparatively minor contributions. DD-GPCE exploits this low-order interaction structure by retaining only basis functions involving at most \(S\) input variables~\cite{lee2023high}. This
truncation reduces the number of basis terms while preserving the interactions most relevant to the response. 

%
Let $m$ be the maximum total polynomial degree and $S$, with $1\le S \le \min(N,m)$, be the maximum retained interaction order.
The corresponding reduced multi-index set is defined as
\begin{equation}
    \mathcal{J}_{S,m}
    :=
    \left\{
        \mathbf{j}\in\mathbb{N}_{0}^{N}
        :
        |\mathbf{j}|\leq m,\;
        \|\mathbf{j}\|_{0}\leq S
    \right\}.
\end{equation}
where $|\mathbf{j}|=\sum_{i=1}^{N}j_i$ denotes the total degree and
$\|\mathbf{j}\|_{0}$ denotes the number of nonzero components of
$\mathbf{j}$. Thus, $\mathcal{J}_{S,m}$ retains polynomial terms involving
at most $S$ input variables. The cardinality of this set is
\begin{equation}
    L_{N,S,m} = 1 + \sum_{s=1}^{S} \binom{N}{s} \binom{m}{s}. \label{eq:cardinality}
\end{equation}
The $S$-variate, $m$th-order DD-GPCE approximation is 
\begin{equation}
    \tilde{y}_{S,m}(\mathbf{X}) = \sum_{i=1}^{L_{N,S,m}} c_i(\mathbf{d}) \Psi_i(\mathbf{X};\mathbf{d})\simeq y(\mathbf{X}), \label{eq:dd_gpce}
\end{equation}
where $\Psi_{i}(\mathbf{X};\mathbf{d})$ is the $i$th multivariate polynomial orthonormal with respect to $f_{\mathbf{X}}(\mathbf{x};\mathbf{d})\mathrm{d}\mathbf{x}$. The corresponding projection coefficient is $c_i(\mathbf{d})=\mathbb{E}_{\mathbf{d}}[y(\mathbf{X})\Psi_{i}(\mathbf{X};\mathbf{d})]$
. The total number of basis functions, $L_{N,S,m}$, is governed by the interaction parameter $S$. For instance, a univariate ($S=1$) DD-GPCE approximation solely accounts for the mean and the independent main effects of each input variable, thereby demanding only $Nm+1$ basis functions. A bivariate ($S=2$) approximation additionally captures the pairwise coupled effects. Consequently, by retaining a low-order interaction level ($S \ll N$), DD-GPCE decomposes the high-dimensional response into a hierarchical structure while preventing the exponential increase of the basis terms.

The expansion coefficients generally involve high-dimensional
integrals. When these integrals cannot be evaluated analytically, the
coefficients can be estimated from input--output data using standard
least squares (SLS). Given $L'$ input samples
$\{\mathbf{x}^{(\ell)}\}_{\ell=1}^{L'}$, define the common design
matrix $\mathbf{A}(\mathbf{d})$ and the response vector
$\mathbf{b}$ for the output as
\begin{equation}
    \mathbf{A}(\mathbf{d}) = \begin{bmatrix} \Psi_1(\mathbf{x}^{(1)};\mathbf{d}) & \dots & \Psi_{L_{N,S,m}}(\mathbf{x}^{(1)};\mathbf{d}) \\ \vdots & \ddots & \vdots \\ \Psi_1(\mathbf{x}^{(L')};\mathbf{d}) & \dots & \Psi_{L_{N,S,m}}(\mathbf{x}^{(L')};\mathbf{d}) \end{bmatrix}, \quad \mathbf{b} = \begin{bmatrix} y(\mathbf{x}^{(1)}) \\ \vdots \\ y(\mathbf{x}^{(L')}) \end{bmatrix}. \label{eq:sls_matrix}
\end{equation}
The coefficient vector $\widehat{\mathbf{c}}(\mathbf{d})=(\widehat{c}_1(\mathbf{d}),\ldots,\widehat{c}_{L_{N,S,m}}(\mathbf{d}))^\top$ minimizes the sum of squared residuals:
\begin{equation}
   \widehat{\mathbf{c}}(\mathbf{d}) = \argmin_{\mathbf{c} \in \mathbb{R}^{L_{N,S,m}}} \| \mathbf{b} - \mathbf{A}(\mathbf{d})\mathbf{c} \|_2^2. \label{eq:sls_obj}
\end{equation}
Provided that the information matrix $\mathbf{A}(\mathbf{d})^\top \mathbf{A}(\mathbf{d})$ is positive-definite, the solution to this problem is analytically computed as
\begin{equation}
    \widehat{\mathbf{c}}(\mathbf{d}) = (\mathbf{A}(\mathbf{d})^\top \mathbf{A}(\mathbf{d}))^{-1} \mathbf{A}(\mathbf{d})^\top \mathbf{b}. \label{eq:sls_sol}
\end{equation}

SLS requires $L'>L_{N,S,m}$. Even with the reduced DD-GPCE basis, globally sampled training sets contain few extreme-tail points, and obtaining them solely from HF simulations is costly. Prior bi-fidelity studies used DD-GPCE to reduce HF evaluations for CVaR estimation under dependent inputs~\cite{LeeKramer2023SMO}. Here, DD-GPCE instead provides the global LF response within an iterative optimization, and selected HF evaluations correct only the surrogate error that affects the tail-risk estimate.

\subsection{Multifidelity adaptation by additive updates}
\label{subsec:multifidelity}

\begin{figure}
    \centering
    \includegraphics[width=1.0\linewidth]{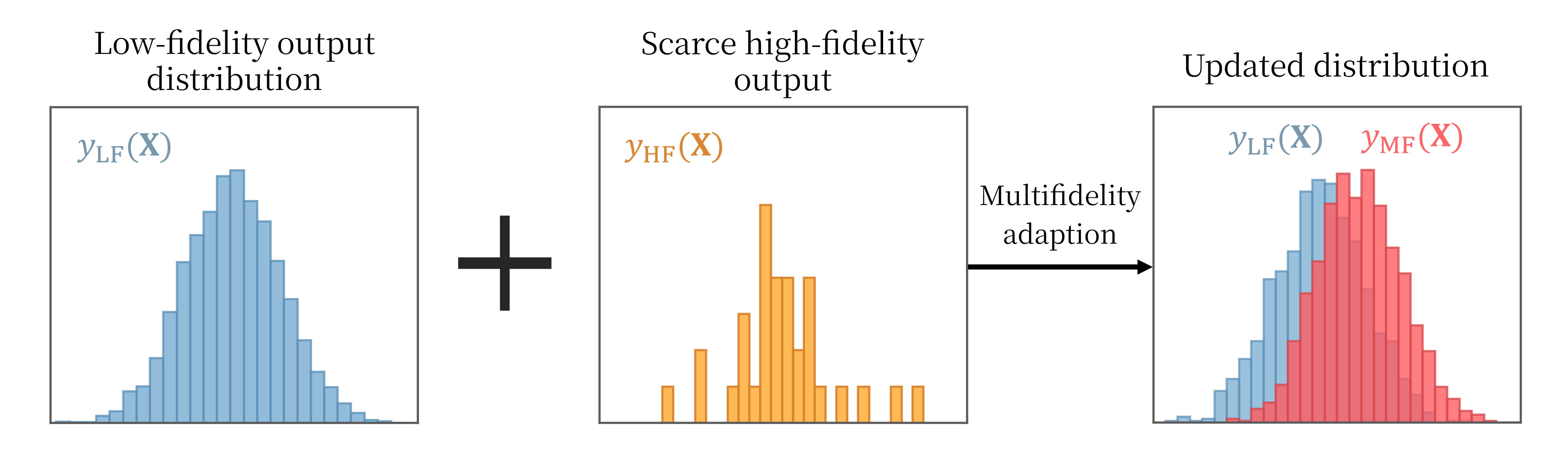}
    \caption{Conceptual workflow of the multifidelity adaptation process: combining a low-fidelity output distribution $f_{\mathrm{LF}}(\mathbf{X})$ with scarce high-fidelity data $f_{\mathrm{HF}}(\mathbf{X})$ to construct an updated and calibrated multifidelity distribution $f_{\mathrm{MF}}(\mathbf{X})$.}
    \label{fig:multifidelity_modeling}
\end{figure}

MF modeling combines an inexpensive LF approximation with limited HF data to correct dominant prediction errors at lower cost~\cite{peherstorfer2018survey}.
A common MF construction represents the HF response as a scaled LF prediction plus an additive discrepancy~\cite{kennedy2000predicting,perdikaris2015multi}. In this study, the scaling coefficient is fixed at unity, yielding the additive model
\begin{equation}
    \hat{y}_{\mathrm{MF}}(\mathbf{x})
    =
    \hat{y}_{\mathrm{LF}}(\mathbf{x})
    +
    \hat{\delta}(\mathbf{x}),
    \label{eq:mf_additive}
\end{equation}
where $\widehat{\delta}(\mathbf{x})$ approximates the LF--HF discrepancy
\begin{equation}
    \delta(\mathbf{x})
    =
    y_{\mathrm{HF}}(\mathbf{x})
    -
    \hat{y}_{\mathrm{LF}}(\mathbf{x}).
\end{equation}
Figure~\ref{fig:multifidelity_modeling} illustrates that  $\hat{y}_{\mathrm{LF}}(\mathbf{x})$ represents the baseline response trend, whereas $\hat{\delta}(\mathbf{x})$ corrects the residuals identified from a limited set of paired LF and HF evaluations. When the LF response and the residual are represented using compatible polynomial bases, the two terms can be combined within a unified multifidelity PCE representation~\cite{ng2012multifidelity,bryson2017all}. 
In this study, we adopt the baseline LF predictor as the DD-GPCE model $\tilde{y}_{S,m}(\mathbf{x})$, such that 
\[
     \widehat{y}_{\mathrm{LF}}(\mathbf{x})
    :=
    \widetilde{y}_{S,m}(\mathbf{x}).
\]

DD-GPCE minimizes average squared prediction error under $f_{\mathbf{X}}(\mathbf{x};\mathbf{d})\mathrm{d}\mathbf{x}$, so high-probability regions dominate the fit. The risk region $\mathcal{G}_{\beta;\mathbf{d}}$ occupies only the upper $100(1-\beta)\%$ of the response distribution; global accuracy therefore does not ensure an accurate CVaR estimate~\cite{Heinkenschloss2020,Jakeman2022,LeeKramer2023}.

Rather than fit a global residual, the proposed method selects HF samples from the current estimate of $\mathcal{G}_{\beta;\mathbf{d}}$ and updates the residual as this region changes during optimization. Figure~\ref{fig:multifidelity_modeling} summarizes how paired LF--HF evaluations correct the DD-GPCE response used for tail-risk estimation.

\begin{problem}
Develop a tail-focused additive MF model that uses scarce HF evaluations to correct DD-GPCE errors directly affecting CVaR, while preserving the global LF structure and avoiding unnecessary HF sampling outside the risk region.
\end{problem}

\section{Risk-averse design optimization via multifidelity tail correction} \label{sec:mf_for_RADO}

This section develops the MF tail-correction method. DD-GPCE represents the global LF response, a confidence-interval-based $\epsilon$-risk region identifies samples that may influence CVaR, and two-stage HF sampling supports a regularized residual correction. The corrected response then provides sampling-based CVaR constraints within a generation-wise optimization loop.

\subsection{Confidence interval-based \texorpdfstring{$\epsilon$}{epsilon}-risk region}\label{subsec:riskregion}

Aleatoric and epistemic uncertainties have distinct roles~\cite{der2009aleatory} in the proposed method. Aleatoric uncertainty represents the inherent variability of the input random vector $\mathbf{X}$ under $\mathrm{P}_{\mathbf{d}}$. This variability induces the response distribution $y(\mathbf{X})$ from which VaR and CVaR are evaluated and cannot be reduced through additional model evaluations. In contrast, epistemic uncertainty arises here from estimating the DD-GPCE coefficients using a finite training data set. 
The confidence interval introduced in this subsection quantifies this coefficient-estimation uncertainty, rather than aleatoric variability in the physical inputs. 

Although DD-GPCE produces a deterministic response at a fixed input $\mathbf{x}$, its estimated coefficients and predictions depend on the available training data. At a fixed design vector $\mathbf{d}$, the SLS model is represented as $\mathbf{b}=\mathbf{A}(\mathbf{d})\mathbf{c}(\mathbf{d})+\mathbf{e}(\mathbf{d})$, $\mathbf{e}(\mathbf{d})\sim\mathcal{N}(\mathbf{0},\sigma^2(\mathbf{d})\mathbf{I})$. The error vector $\mathbf{e}(\mathbf{d})$ is a working statistical representation of unresolved regression and truncation errors; the error does not represent the aleatoric variability of $\mathbf{X}$. Under the classical assumption of independent, homoscedastic Gaussian errors and a full-column-rank design matrix, this model provides the standard coefficient covariance, unbiased residual-variance estimator, and finite-sample Student's $t$ confidence interval~\cite{sacks1989design,seber2003linear,lange1989robust}. 

The SLS coefficient estimator and its estimation error are $\widehat{\mathbf{c}}(\mathbf{d}) = (\mathbf{A}(\mathbf{d})^\top \mathbf{A}(\mathbf{d}))^{-1} \mathbf{A}(\mathbf{d})^\top \mathbf{b}$ and $\widehat{\mathbf{c}}(\mathbf{d}) - \mathbf{c}(\mathbf{d}) = (\mathbf{A}(\mathbf{d})^\top \mathbf{A}(\mathbf{d}))^{-1} \mathbf{A}(\mathbf{d})^\top \mathbf{e}(\mathbf{d})$. Consequently, the coefficient covariance is
\begin{equation}
    \operatorname{Var}[\widehat{\mathbf{c}}(\mathbf{d})] = \mathbb{E}[(\widehat{\mathbf{c}}(\mathbf{d}) - \mathbf{c}(\mathbf{d}))(\widehat{\mathbf{c}}(\mathbf{d}) - \mathbf{c}(\mathbf{d}))^\top] = \sigma^2(\mathbf{d}) (\mathbf{A}(\mathbf{d})^\top \mathbf{A}(\mathbf{d}))^{-1}. \label{eq:coeff_variance}
\end{equation}

To estimate the residual variance $\sigma^2(\mathbf{d})$, we use the residual sum of squares. The hat matrix is $\mathbf{H}(\mathbf{d}) = \mathbf{A}(\mathbf{d})(\mathbf{A}(\mathbf{d})^\top \mathbf{A}(\mathbf{d}))^{-1} \mathbf{A}(\mathbf{d})^\top$, and the observed residual vector is $\boldsymbol{\delta}(\mathbf{d}) = (\mathbf{I} - \mathbf{H}(\mathbf{d}))\mathbf{b} = (\mathbf{I} - \mathbf{H}(\mathbf{d}))\mathbf{e}(\mathbf{d})$. The expectation of the residual sum of squares is
\begin{equation}
    \mathbb{E}[\boldsymbol{\delta}(\mathbf{d})^\top \boldsymbol{\delta}(\mathbf{d})] = \mathbb{E}[\mathbf{e}(\mathbf{d})^\top (\mathbf{I} - \mathbf{H}(\mathbf{d})) \mathbf{e}(\mathbf{d})] = \sigma^2(\mathbf{d}) \text{tr}(\mathbf{I} - \mathbf{H}(\mathbf{d})) = \sigma^2(\mathbf{d}) (L' - L_{N,S,m}). \label{eq:expected_rss}
\end{equation}
Therefore, an unbiased estimator for the residual variance, $\hat{\sigma}^2(\mathbf{d})$, is defined as
\begin{equation}
    \hat{\sigma}^2(\mathbf{d}) = \frac{\|\mathbf{b} - \mathbf{A}(\mathbf{d})\hat{\mathbf{c}}(\mathbf{d})\|^2}{L' - L_{N,S,m}}, \label{eq:unbiased_variance}
\end{equation}
where $L' - L_{N,S,m}$ represents the degrees of freedom corresponding to the sample size $L'$ and the number of basis functions $L_{N,S,m}$.

For a generic input $\mathbf{x}$, the estimated prediction variance $\hat{s}^2(\mathbf{x};\mathbf{d})$ originating from the uncertainty of the coefficient estimation is evaluated as
\begin{equation}
    \hat{s}^2(\mathbf{x};\mathbf{d}) = \hat{\sigma}^2(\mathbf{d}) \mathbf{\Psi}^\top(\mathbf{x};\mathbf{d}) (\mathbf{A}(\mathbf{d})^\top \mathbf{A}(\mathbf{d}))^{-1} \mathbf{\Psi}(\mathbf{x};\mathbf{d}). \label{eq:predictive_variance}
\end{equation}
Using this variance, the $(1-\alpha)100\%$ confidence interval $\epsilon_{S,m}^{\alpha}(\mathbf{x})$ associated with the t-distribution is defined as
\begin{equation}
    \epsilon_{S,m}^{\alpha}(\mathbf{x};\mathbf{d}) = t_{\nu, 1-\alpha/2} \hat{s}(\mathbf{x};\mathbf{d}), \label{eq:confidence_width}
\end{equation}
where $t_{\nu, 1-\alpha/2}$ is the critical value of the t-distribution with $\nu = L' - L_{N,S,m}$ degrees of freedom.

We define the $(1-\alpha)100\%$ CI-based $\epsilon$-risk region as a discrete candidate set. To reduce the chance of excluding relevant tail samples, we estimate $\widehat{\mathrm{VaR}}_{\beta;\mathbf{d}}$ from the lower confidence bounds and retain samples whose upper confidence bounds exceed that threshold:
\begin{equation}
\begin{aligned}
\widehat{\mathcal{G}}_{\beta;\mathbf{d}}^{\alpha}:=\Bigl\{
\mathbf{x}^{(j)} \in \mathbb{A}^N \;\Big|\;& \tilde{y}_{S,m}(\mathbf{x}^{(j)})
+\epsilon_{S,m}^{\alpha}(\mathbf{x}^{(j)};\mathbf{d}) \\ &\ge
\widehat{\mathrm{VaR}}_{\beta;\mathbf{d}} \left[\tilde{y}_{S,m}(\mathbf{X})
- \epsilon_{S,m}^{\alpha}(\mathbf{X};\mathbf{d})\right],\quad j=1,\ldots,L \Bigr\},
&\bar{L}:=\left|\widehat{\mathcal{G}}_{\beta;\mathbf{d}}^{\alpha}\right|<L.
\end{aligned}
\label{eq:risk_region_discrete}
\end{equation}
Here, $L$ is the number of candidate samples and $\bar{L}$ is the number retained. The construction expands the estimated tail according to surrogate uncertainty, and subsequent HF selection is restricted to $\widehat{\mathcal{G}}_{\beta;\mathbf{d}}^{\alpha}$.

\vspace{\baselineskip}
\noindent \textbf{Remark 1.} The probabilistic CI does not guarantee that $\widehat{\mathcal{G}}_{\beta;\mathbf{d}}^{\alpha}$ always contains the true risk region~\cite{Heinkenschloss2020}; it provides a conservative candidate set that reduces the chance of missing influential tail responses.
\begin{algorithm}[htb]
\caption{Construction of the CI-based $\epsilon$-risk region using DD-GPCE}
\label{alg:risk_region}
\begin{algorithmic}[1]
\renewcommand{\algorithmicrequire}{\textbf{Input:}}
\renewcommand{\algorithmicensure}{\textbf{Output:}}
\Require Fixed design vector $\mathbf{d}$; candidate samples $\mathbf{x}^{(l)} = (x_1^{(l)}, \dots, x_N^{(l)})^\top$, $l = 1, \dots, L$, generated from $f_{\mathbf{X}}(\mathbf{x};\mathbf{d})\mathrm{d}\mathbf{x}$ via MCS, quasi MCS, or Latin hypercube sampling; normalized sample weights 
$\{p^{(l)}\}_{l=1}^{L}$ satisfying
$p^{(l)}\geq0$ and $\sum_{l=1}^{L}p^{(l)}=1$; training input-output data set $\{\mathbf{x}^{(q)}, y(\mathbf{x}^{(q)})\}_{q=1}^{L'}$, $L' < L$;
confidence parameter $\alpha\in(0,1)$;
risk level $\beta\in(0,1)$.
\Ensure CI-based $\epsilon$-risk region $\widehat{\mathcal{G}}_{\beta;\mathbf{d}}^{\alpha}$.

\Procedure{Calculate DD-GPCE}{$\tilde{y}_{S,m}(\mathbf{x})$}
    \State Create the measure-consistent orthonormal basis vector $\mathbf{\Psi}_{S,m}(\mathbf{x};\mathbf{d})$ as described in Supplementary Section~S2.
    \State Compute the DD-GPCE coefficient vector $\widehat{\mathbf{c}}(\mathbf{d})$ using Eq.~\eqref{eq:sls_sol}.
\EndProcedure

\State Evaluate $\tilde{y}_{S,m}(\mathbf{x}^{(l)})$, $\hat{s}^2(\mathbf{x}^{(l)};\mathbf{d})$, and the CI half-width $\epsilon_{S,m}^{\alpha}(\mathbf{x}^{(l)};\mathbf{d})$ for $l=1,\ldots,L$ using Eq.~\eqref{eq:confidence_width}.

\State Sort the lower confidence bounds in descending order and relabel the samples accordingly:
\[
    \tilde{y}_{S,m}(\mathbf{x}^{(1)}) - \epsilon_{S,m}^{\alpha}(\mathbf{x}^{(1)};\mathbf{d}) \ge \cdots \ge \tilde{y}_{S,m}(\mathbf{x}^{(L)}) - \epsilon_{S,m}^{\alpha}(\mathbf{x}^{(L)};\mathbf{d}).
\]

\State Compute an index $\bar{k}_\beta$ such that
\[
    \sum_{l=1}^{\bar{k}_\beta - 1} p^{(l)} \le 1 - \beta < \sum_{l=1}^{\bar{k}_\beta} p^{(l)}.
\]

\State Set
\[
    \widehat{\mathrm{VaR}}_{\beta;\mathbf{d}}[\tilde{y}_{S,m}(\mathbf{X}) - \epsilon_{S,m}^{\alpha}(\mathbf{X};\mathbf{d
    })] = \tilde{y}_{S,m}(\mathbf{x}^{(\bar{k}_\beta)})-\epsilon_{S,m}^{\alpha}(\mathbf{x}^{(\bar{k}_\beta)};\mathbf{d}).
\]

\State Define the risk-region index set:
\[
    \mathcal{I}_{\beta;\mathbf{d}}^{\alpha}=\left\{j\in\{1,\dots,L\}\;\middle|\;
    \tilde{y}_{S,m}(\mathbf{x}^{(j)}) + \epsilon_{S,m}^{\alpha}(\mathbf{x}^{(j)};\mathbf{d}) \geq\widehat{\mathrm{VaR}}_{\beta;\mathbf{d}}[\tilde{y}_{S,m}(\mathbf{X}) - \epsilon_{S,m}^{\alpha}(\mathbf{X};\mathbf{d})]\right\}.
\]

\State Determine the $(1-\alpha)100\%$ CI-based $\epsilon$-risk region by a discrete set:
\begin{align*}
    \widehat{\mathcal{G}}_{\beta;\mathbf{d}}^{\alpha}
    = \left\{\mathbf{x}^{(j)}\;\middle|\;j\in\mathcal{I}_{\beta;\mathbf{d}}^{\alpha}
    \right\}, \qquad \bar L=\left|\mathcal{I}_{\beta;\mathbf{d}}^{\alpha}\right|.
\end{align*}
\end{algorithmic}
\end{algorithm}

\subsection{Infill criterion for tail correction} \label{subsec:criterion}

The proposed two-stage adaptive sampling strategy allocates the HF budget sequentially within $\widehat{\mathcal{G}}_{\beta;\mathbf{d}}^\alpha$, balancing exploration and exploitation~\cite{crombecq2011novel}. Stage~1 targets high DD-GPCE uncertainty; Stage~2 uses the intermediate correction to locate predictions that may need further adjustment.

In Stage~1, we evaluate the $N_1$ candidates with the largest prediction variances $\widehat{s}^{2}(\mathbf{x};\mathbf{d})$ from Eq.~\eqref{eq:predictive_variance}. Their LF--HF residuals define the intermediate response $\widehat{y}_{\mathrm{MF}}^{(1)}$ through the additive correction in Section~\ref{subsec:correction}.

The second stage shifts the emphasis toward correction-informed exploitation while retaining an uncertainty-based weighting. The acquisition score is
\begin{equation}
    a_{\mathrm{shift}}(\mathbf{x};\mathbf{d}) = \left| \widehat{y}_{\mathrm{MF}}^{(1)}(\mathbf{x}) - \widehat{y}_{S,m}(\mathbf{x}) \right| \times \widehat{s}(\mathbf{x};\mathbf{d}), \label{eq:shift_metric}
\end{equation}
The absolute shift measures the change induced by the Stage~1 correction, while multiplication by $\widehat{s}(\mathbf{x};\mathbf{d})$ favors locations that also remain weakly supported by the original DD-GPCE data.

Stage~2 selects the $N_2$ unevaluated candidates with the largest $a_{\mathrm{shift}}(\mathbf{x};\mathbf{d})$. Combining both stages yields the final correction from $N_1+N_2$ HF evaluations: $N_1$ uncertainty-driven samples and $N_2$ correction-informed samples.

\subsection{Residual-based tail correction via DD-GPCE multifidelity}
\label{subsec:correction}

Let $\mathcal{I}_{\beta;\mathbf{d}}^{\alpha}=\{i\in\{1,\ldots,L\}\mid\mathbf{x}^{(i)}\in\widehat{\mathcal{G}}_{\beta;\mathbf{d}}^{\alpha}\}$. Section~\ref{subsec:criterion} selects ordered indices $i_1<\cdots<i_{n_c}$ from this set, with $n_c=N_1$ for the intermediate correction and $n_c=N_1+N_2$ for the final correction.

MF PCE represents the HF response as the sum of a baseline LF polynomial approximation and an additive LF--HF discrepancy~\cite{palar2016multi}. In this method,
\begin{equation}
    y_{\mathrm{HF}}(\mathbf{x})=\tilde{y}_{S,m}(\mathbf{x})+\delta(\mathbf{x};\mathbf{d}).
    \label{eq:exact_response}
\end{equation}
where $y_{\mathrm{HF}}(\mathbf{x})\equiv y(\mathbf{x})$ is the HF response, $\tilde{y}_{S,m}(\mathbf{x})$ is the baseline DD-GPCE prediction, and $\delta(\mathbf{x};\mathbf{d}):=y_{\mathrm{HF}}(\mathbf{x})-\widetilde{y}_{S,m}(\mathbf{x})$ is the LF--HF discrepancy.

The selected HF input locations and their corresponding responses are collected in
\begin{equation}
    \mathbf{X}_{c}=
    \begin{bmatrix}
        \mathbf{x}^{(i_1)} &  \cdots & \mathbf{x}^{(i_{n_c})}
    \end{bmatrix}^{\top}\in\mathbb{R}^{n_c\times N},
    \label{eq:selected_input_matrix}
\end{equation}
and 
\begin{equation}
    \mathbf{y}_{c}=
    \begin{bmatrix}
        y(\mathbf{x}^{(i_1)}) & \cdots &
        y(\mathbf{x}^{(i_{n_c})})
    \end{bmatrix}^{\top}
    \in \mathbb{R}^{n_c},
    \label{eq:selected_exact_vector}
\end{equation}
respectively. The ordered residual vector is then given by
\begin{equation}
    \boldsymbol{\delta}_{c}=\mathbf{y}_{c}-
    \begin{bmatrix}
        \tilde{y}_{S,m}(\mathbf{x}^{(i_1)}) &
    \cdots &
        \tilde{y}_{S,m}(\mathbf{x}^{(i_{n_c})})
    \end{bmatrix}^{\top}.
    \label{eq:exact_residual}
\end{equation}
We approximate the discrepancy $\delta(\mathbf{x};\mathbf{d})$ using the same measure-consistent orthonormal basis as the baseline DD-GPCE:
\begin{equation}
    \tilde{\delta}_{S,m}(\mathbf{x};\mathbf{d})=\boldsymbol{\Psi}_{S,m}(\mathbf{x};\mathbf{d})^\top \boldsymbol{\gamma}(\mathbf{d}),
    \label{eq:residual_expansion}
\end{equation}
where
\begin{equation}
    \boldsymbol{\Psi}_{S,m}(\mathbf{x};\mathbf{d})
    :=
    \begin{bmatrix}
        \Psi_1(\mathbf{x};\mathbf{d}) &
        \cdots &
        \Psi_{L_{N,S,m}}(\mathbf{x};\mathbf{d})
    \end{bmatrix}^{\top}
    \label{eq:residual_basis_vector}
\end{equation}
is the DD-GPCE basis vector and $\boldsymbol{\gamma}(\mathbf{d})\in\mathbb{R}^{L_{N,S,m}}$ is the correction coefficient vector.
At the selected tail samples, the discrepancy approximation satisfies
\begin{equation}
    \boldsymbol{\delta}_{c}(\mathbf{d})
    \approx
    \mathbf{\Psi}_{c;S,m}(\mathbf{d})\boldsymbol{\gamma}(\mathbf{d}),
    \label{eq:residual_approx}
\end{equation}
where the ordered residual design matrix $\mathbf{\Psi}_{c;S,m}(\mathbf{d})\in\mathbb{R}^{n_c\times L_{N,S,m}}$ is defined by
\begin{equation}
    \mathbf{\Psi}_{c;S,m}(\mathbf{d})=
    \begin{bmatrix}
        \boldsymbol{\Psi}_{S,m}(\mathbf{x}^{(i_1)};\mathbf{d})^\top \\
 \vdots \\
        \boldsymbol{\Psi}_{S,m}(\mathbf{x}^{(i_{n_c})};\mathbf{d})^\top
    \end{bmatrix}.
    \label{eq:residual_design_matrix}
\end{equation}
Consequently, the exact response is approximated by the tail-corrected DD-GPCE model $\tilde{y}_{\mathrm{MF}}(\mathbf{x})$, formulated as an additive update to the baseline LF prediction:
\begin{equation}
    \tilde{y}_{\mathrm{MF}}(\mathbf{x})=\tilde{y}_{S,m}(\mathbf{x})
    +\tilde{\delta}_{S,m}(\mathbf{x};\mathbf{d}).
    \label{eq:corrected_model}
\end{equation}

Because the number of HF tail evaluations is intentionally limited, $n_c<L_{N,S,m}$ typically holds. Consequently, the SLS regression does not provide a unique coefficient vector. The correction coefficients are therefore estimated using Tikhonov regularization~\cite{tikhonov1963solution}:
\begin{equation}   \boldsymbol{\gamma}_{\lambda}(\mathbf{d})=\arg\min_{\boldsymbol{\gamma}}
    \left\{\left\|\mathbf{\Psi}_{c;S,m}(\mathbf{d})\boldsymbol{\gamma}-
    \boldsymbol{\delta}_{c}(\mathbf{d})\right\|_2^2+\lambda\left\|\boldsymbol{\gamma}
    \right\|_2^2\right\},\qquad\lambda>0,
    \label{eq:tikhonov_obj}
\end{equation}
where $\lambda$ is the regularization parameter and the identity matrix is adopted as the Tikhonov matrix, i.e., $\mathbf{D}=\mathbf{I}$. The corresponding closed-form solution is
\begin{equation}
    \boldsymbol{\gamma}_{\lambda}(\mathbf{d})=\left(\mathbf{\Psi}_{c;S,m}^{\top}(\mathbf{d})
    \mathbf{\Psi}_{c;S,m}(\mathbf{d})+\lambda\mathbf{I}\right)^{-1}\mathbf{\Psi}_{c;S,m}^{\top}(\mathbf{d})
    \boldsymbol{\delta}_{c}(\mathbf{d}).
    \label{eq:tikhonov_sol}
\end{equation}

For any $\lambda>0$, the matrix
$\mathbf{\Psi}_{c;S,m}^{\top}(\mathbf{d})\mathbf{\Psi}_{c;S,m}(\mathbf{d})+\lambda\mathbf{I}$
is positive definite. Therefore, its inverse exists, and the regularized correction coefficients $\boldsymbol{\gamma}_{\lambda}(\mathbf{d})$ are uniquely determined. The regularization term also suppresses excessively large coefficient magnitudes that could otherwise arise from the sparse and underdetermined tail-data system. The selected exact samples therefore serve as local calibration data for adjusting the magnitude and direction of the surrogate discrepancy within the risk region.

Because the baseline and discrepancy expansions use the same
polynomial basis, the corrected response can be written as a single DD-GPCE expansion:
\begin{equation}
    \tilde{y}_{\mathrm{MF}}(\mathbf{x})=\sum_{i=1}^{L_{N,S,m}}
    \left(\widehat{c}_{i}(\mathbf{d})+\gamma_{\lambda,i}(\mathbf{d})\right)\Psi_i(\mathbf{x};\mathbf{d}),
    \label{eq:final_correction}
\end{equation}
where $\widehat{\mathbf{c}}(\mathbf{d})$ is the baseline coefficient vector and $\gamma_{\lambda,i}(\mathbf{d})$ is the $i$th element of $\boldsymbol{\gamma}_{\lambda}(\mathbf{d})$. This unified representation avoids the separate evaluation of multiple surrogate models and streamlines the subsequent sampling-based CVaR assessment.

\subsection{Sampling-based CVaR estimation using tail-corrected DD-GPCE} \label{subsec:estimation}
\begin{algorithm}[htb]
\caption{Sampling-based estimation of $\mathrm{VaR}_{\beta;\mathbf{d}}$ and $\mathrm{CVaR}_{\beta;\mathbf{d}}$ using the tail-corrected DD-GPCE}
\label{alg:cvar_estimation}
\begin{algorithmic}[1]
\renewcommand{\algorithmicrequire}{\textbf{Input:}}
\renewcommand{\algorithmicensure}{\textbf{Output:}}

\Require Fixed design vector $\mathbf{d}$, candidate input samples $\{\mathbf{x}^{(l)}\}_{l=1}^{L}$, where $\mathbf{x}^{(l)}=(x_1^{(l)},\dots,x_N^{(l)})^\top$ and $L\gg1$, generated using MCS, quasi-MCS, or Latin hypercube sampling; corresponding normalized probability weights $\{p^{(l)}\}_{l=1}^{L}$ satisfying $p^{(l)}\geq0$ and $\sum_{l=1}^{L}p^{(l)}=1$; baseline coefficient vector $\widehat{\mathbf{c}}(\mathbf{d})\in\mathbb{R}^{L_{N,S,m}}$; regularized correction coefficient vector $\boldsymbol{\gamma}_{\lambda}(\mathbf{d})\in\mathbb{R}^{L_{N,S,m}}$ obtained from Section~\ref{subsec:correction}; risk level $\beta\in(0,1)$.

\Ensure Estimates $\widehat{\mathrm{VaR}}_{\beta;\mathbf{d}}[\tilde{y}_{\mathrm{MF}}(\mathbf{X})]$ and $\widehat{\mathrm{CVaR}}_{\beta;\mathbf{d}}[\tilde{y}_{\mathrm{MF}}(\mathbf{X})]$.

\For{$l=1,\dots,L$}
    \State Evaluate 
    $\tilde{y}_{\mathrm{MF}}(\mathbf{x}^{(l)})
    \gets
    \sum_{i=1}^{L_{N,S,m}}
    (\widehat{c}_{i}(\mathbf{d})+\gamma_{\lambda,i}(\mathbf{d}))
    \Psi_i(\mathbf{x}^{(l)};\mathbf{d})$.
\EndFor

\State Determine a permutation
\[
    \boldsymbol{\pi}=(\pi_1,\dots,\pi_L)
\]
of $\{1,\dots,L\}$ such that the output samples are arranged in descending order:
\[
    \tilde{y}_{\mathrm{MF}}(\mathbf{x}^{(\pi_1)}) \geq
    \tilde{y}_{\mathrm{MF}}(\mathbf{x}^{(\pi_2)}) \geq
    \cdots \geq \tilde{y}_{\mathrm{MF}}(\mathbf{x}^{(\pi_L)}).
\]

\State Reorder the probability weights according to the same permutation, resulting in
\[
    p^{(\pi_1)},p^{(\pi_2)},\dots,p^{(\pi_L)}.
\]

\State Determine the weighted upper-tail index $k_{\beta}\in\{1,\dots,L\}$ satisfying
\[
    \sum_{l=1}^{k_{\beta}-1}p^{(\pi_l)} \leq 1-\beta <
    \sum_{l=1}^{k_{\beta}}p^{(\pi_l)}.
\]

\State Estimate the Value-at-Risk as
\[
    \widehat{\mathrm{VaR}}_{\beta;\mathbf{d}}[\tilde{y}_{\mathrm{MF}}(\mathbf{X})]
    =\tilde{y}_{\mathrm{MF}}(\mathbf{x}^{(\pi_{k_{\beta}})}).
\]

\State Estimate the Conditional Value-at-Risk as
\[
\begin{aligned}
    \widehat{\mathrm{CVaR}}_{\beta;\mathbf{d}}[\tilde{y}_{\mathrm{MF}}(\mathbf{X})]=
    \widehat{\mathrm{VaR}}_{\beta;\mathbf{d}}[\tilde{y}_{\mathrm{MF}}(\mathbf{X})]+
    \frac{1}{1-\beta}\sum_{l=1}^{L}p^{(l)}\left( \tilde{y}_{\mathrm{MF}}(\mathbf{x}^{(l)}) - \widehat{\text{VaR}}_{\beta;\mathbf{d}} [\tilde{y}_{\mathrm{MF}}(\mathbf{X})] \right)^+.
\end{aligned}
\]
\end{algorithmic}
\end{algorithm}
Only about $(1-\beta)L$ of $L$ equally weighted samples lie in the upper tail, so direct HF evaluation is impractical at high risk levels~\cite{hong2014monte}. We instead evaluate the corrected model at samples $\{\mathbf{x}^{(l)}\}_{l=1}^{L}$ with normalized weights $p^{(l)}\geq0$ and $\sum_l p^{(l)}=1$ ($p^{(l)}=1/L$ for equal weights). Sorting the corrected responses and their weights gives the empirical VaR; the positive-part form then gives CVaR:
\begin{equation}
    \widehat{\mathrm{CVaR}}_{\beta;\mathbf{d}} [\tilde{y}_{\mathrm{MF}}(\mathbf{X})] = \widehat{\mathrm{VaR}}_{\beta;\mathbf{d}} [\tilde{y}_{\mathrm{MF}}(\mathbf{X})] + \frac{1}{1-\beta}\sum_{l=1}^{L}p^{(l)} \left( \tilde{y}_{\mathrm{MF}}(\mathbf{x}^{(l)}) - \widehat{\mathrm{VaR}}_{\beta;\mathbf{d}} [\tilde{y}_{\mathrm{MF}}(\mathbf{X})] \right)^+. \label{eq:cvar_estimate}
\end{equation}
Here, $(\cdot)^+=\max(0,\cdot)$. Algorithm~\ref{alg:cvar_estimation} gives the complete weighted procedure, including the upper-tail index and probability mass at the empirical VaR threshold. After tail correction, this many-query estimate requires no additional HF evaluations~\cite{LeeKramer2023}.
%
\subsection{Integration with the risk-averse optimization loop}
\label{subsec:optimization}

We solve problem~\eqref{eq:multioutput_risk_averse_design} by evaluating the objective and tail-corrected CVaR constraints for each candidate design. At generation $k$, the current surrogate $\widetilde{y}_{\mathrm{MF},j}^{(k)}$ and $\mathbf{X}\sim\mathrm{P}_{\mathbf{d}}$ define
\begin{equation}
    \widehat{g}_{j}^{(k)}(\mathbf{d}) := \widehat{\mathrm{CVaR}}_{\beta;\mathbf{d}}
    \!\left[\widetilde{y}_{\mathrm{MF},j}^{(k)}(\mathbf{X})\right]
    - \tau_j, \qquad j=1,\ldots,J.
    \label{eq:estimated_cvar_constraint}
\end{equation}
The optimizer compares candidates using $c(\mathbf{d})$ and
\[
\widehat{\mathbf{g}}^{(k)}(\mathbf{d})
=\left(\widehat{g}_{1}^{(k)}(\mathbf{d}),\ldots,
\widehat{g}_{J}^{(k)}(\mathbf{d})\right)^{\top},
\]
where $\widehat{g}_{j}^{(k)}(\mathbf{d})\leq0$ indicates estimated feasibility.

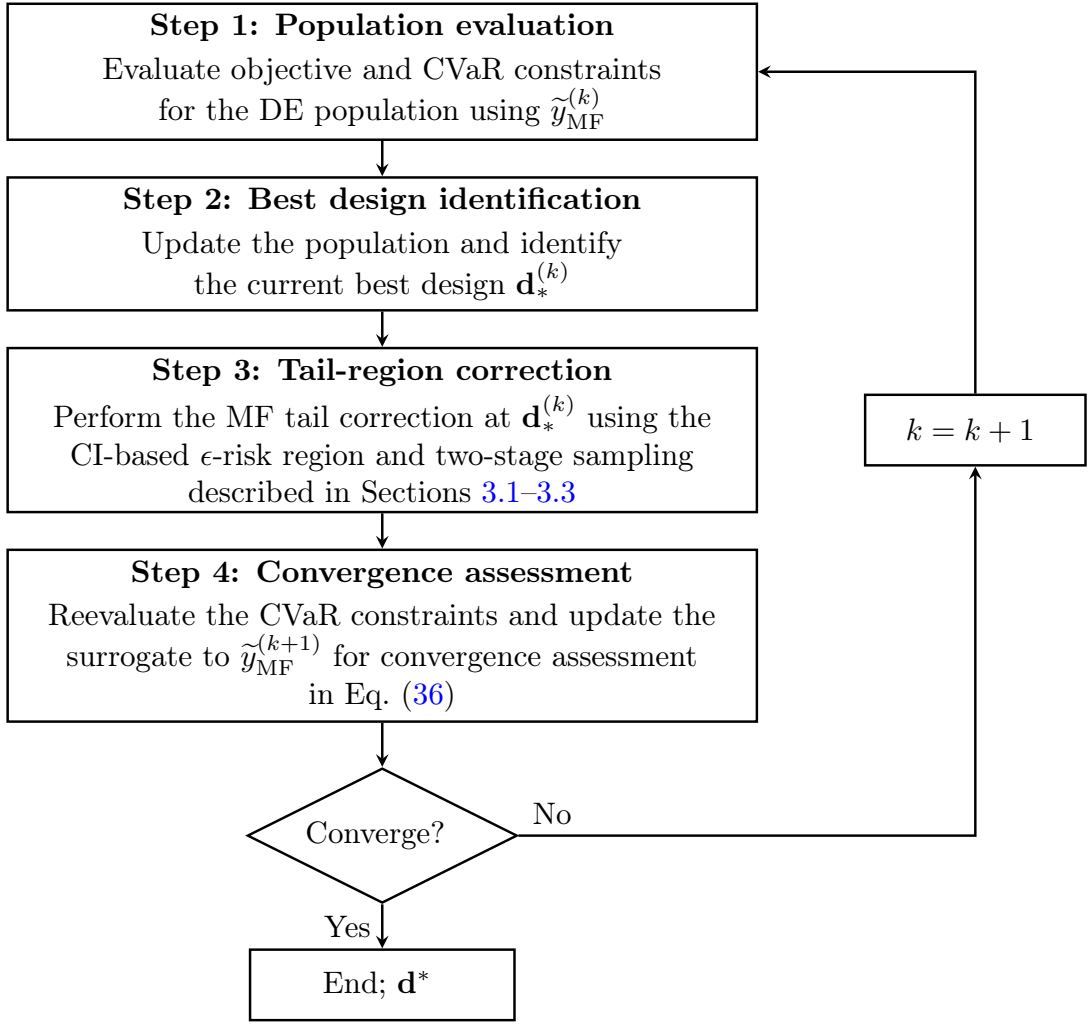
\begin{figure}[htb!]
    \centering
    \resizebox{0.94\textwidth}{!}{
        \begin{tikzpicture}[
            node distance=0.4cm and 1.2cm,
            process/.style={rectangle, draw=black, thick, minimum width=8.5cm, minimum height=0.8cm, align=center, fill=white},
            decision/.style={diamond, draw=black, thick, aspect=2, minimum width=3cm, align=center, fill=white},
            io/.style={rectangle, draw=black, thick, minimum width=3cm, minimum height=0.8cm, align=center, fill=white},
            arrow/.style={thick, ->, >=stealth},
            sideproc/.style={rectangle, draw=black, thick, minimum width=2.5cm, minimum height=0.8cm, align=center, fill=white}
        ]

            \node (s1) [process] {
                \textbf{Step 1: Population evaluation}\\[0.7mm]
                Evaluate objective and CVaR constraints\\for the DE population using $\widetilde{y}_{\mathrm{MF}}^{(k)}$
            };

            \node (s2) [process, below=of s1] {
                \textbf{Step 2: Best design identification}\\[0.7mm]
                Update the population and identify\\the current best design $\mathbf{d}_{*}^{(k)}$
            };

            \node (s3) [process, below=of s2] {
                \textbf{Step 3: Tail-region correction}\\[0.7mm]
                Perform the MF tail correction at $\mathbf{d}_{*}^{(k)}$ using the\\CI-based $\epsilon$-risk region and two-stage sampling\\described in Sections~\ref{subsec:riskregion}--\ref{subsec:correction}
            };

            \node (s4) [process, below=of s3] {
                \textbf{Step 4: Convergence assessment}\\[0.7mm]
                Reevaluate the CVaR constraints and update the\\surrogate to $\widetilde{y}_{\mathrm{MF}}^{(k+1)}$ for convergence assessment\\in Eq.~\eqref{eq:optimization_convergence}
            };

            \node (s5) [decision, below=0.5cm of s4] {
                Converge?
            };

            \node (end) [io, below=0.5cm of s5] {
                End; $\mathbf{d}^*$
            };

            \node (update) [sideproc, right=1.2cm of s3] {
                $k = k + 1$
            };

            \draw [arrow] (s1) -- (s2);
            \draw [arrow] (s2) -- (s3);
            \draw [arrow] (s3) -- (s4);
            \draw [arrow] (s4) -- (s5);
            \draw [arrow] (s5) -- node[anchor=east] {Yes} (end);

            \draw [arrow] (s5.east) -| node[anchor=south, xshift=-4.8cm] {No} (update.south);
            \draw [arrow] (update.north) |- (s1.east);

        \end{tikzpicture}
    } 
    \vspace{10pt}
    \caption{Generation-wise interaction between DE and MF tail correction.}
    \label{fig:de_mf_tail_correction_flow}
\end{figure}

We use Differential Evolution (DE), a derivative-free population-based optimizer~\cite{storn1997differential}, because the sampling, risk-region identification, and residual correction do not yield readily available analytical CVaR gradients and the responses may be nonlinear and nonconvex. Other derivative-free global optimizers can use the same estimator. Figure~\ref{fig:de_mf_tail_correction_flow} shows the generation-wise coupling.

Initial HF data train one global DD-GPCE over the design and uncertainty domain. At the start of generation $k$, $\widetilde{y}_{\mathrm{MF},j}^{(k)}$ denotes its current tail-corrected form. We reuse fixed underlying random realizations, transform them according to $\mathrm{P}_{\mathbf{d}}$ for each candidate, and estimate the response distribution and CVaR constraints. We keep the surrogate fixed while evaluating all trial vectors in a generation and apply the correction once before starting the next generation.

Within each generation, SciPy's \texttt{best1bin} strategy applies mutation, crossover, selection, and immediate population updating~\cite{virtanen2020scipy}. After all trial vectors have been evaluated, we select $\mathbf{d}_{*}^{(k)}$, transform the common random realizations under $\mathrm{P}_{\mathbf{d}_{*}^{(k)}}$, construct its CI-based $\epsilon$-risk region, and perform the two-stage HF correction. The acquired HF samples are appended to the global training set to obtain $\widetilde{y}_{\mathrm{MF},j}^{(k+1)}$. This generation-wise management follows established surrogate-assisted evolutionary practice~\cite{wang2018global}; nonlinear risk constraints use Lampinen's rule as implemented in DE~\cite{lampinen2002constraint}.

After each update, we reevaluate the current best design. Starting at $k_{\min}=10$, let $\widehat{\mathrm{CVaR}}_{j}^{(k,-)}$ and $\widehat{\mathrm{CVaR}}_{j}^{(k,+)}$ be the $j$th estimates before and after correction, and let $\overline{\mathbf{d}}_{*}^{(k)}$ be the normalized best design. We declare convergence when
\begin{equation}
\begin{gathered}
\left\|\overline{\mathbf{d}}_{*}^{(k)}
-\overline{\mathbf{d}}_{*}^{(k-1)}\right\|_{\infty}
< \epsilon_{\mathbf d},
\\[3mm]
\max_{j=1,\ldots,J} \frac{\left|\widehat{\mathrm{CVaR}}_{j}^{(k,-)}-\widehat{\mathrm{CVaR}}_{j}^{(k-1,+)}\right|}{|\tau_j|}
< \epsilon_{\mathrm{CVaR}},
\qquad
\max_{j=1,\ldots,J} \frac{\left|\widehat{\mathrm{CVaR}}_{j}^{(k,+)}-\widehat{\mathrm{CVaR}}_{j}^{(k,-)}\right|}{|\tau_j|}
< \epsilon_{\mathrm{corr}}.
\end{gathered}
\label{eq:optimization_convergence}
\end{equation}
The three conditions bound the design change, the intergeneration CVaR change, and the correction-induced CVaR change; all constraints must also remain satisfied. We set $\epsilon_{\mathbf d}=10^{-3}$ and $\epsilon_{\mathrm{CVaR}}=\epsilon_{\mathrm{corr}}=10^{-3}$, corresponding to $0.1\%$ of $|\tau_j|$ for the CVaR tests. Optimization stops when either DE's criterion ($\mathrm{tol}=0.01$) or Eq.~\eqref{eq:optimization_convergence} is satisfied and returns the best feasible design $\mathbf{d}^{\ast}$.

\subsection{Overall risk-averse optimization algorithm} \label{subsec:algorithm}
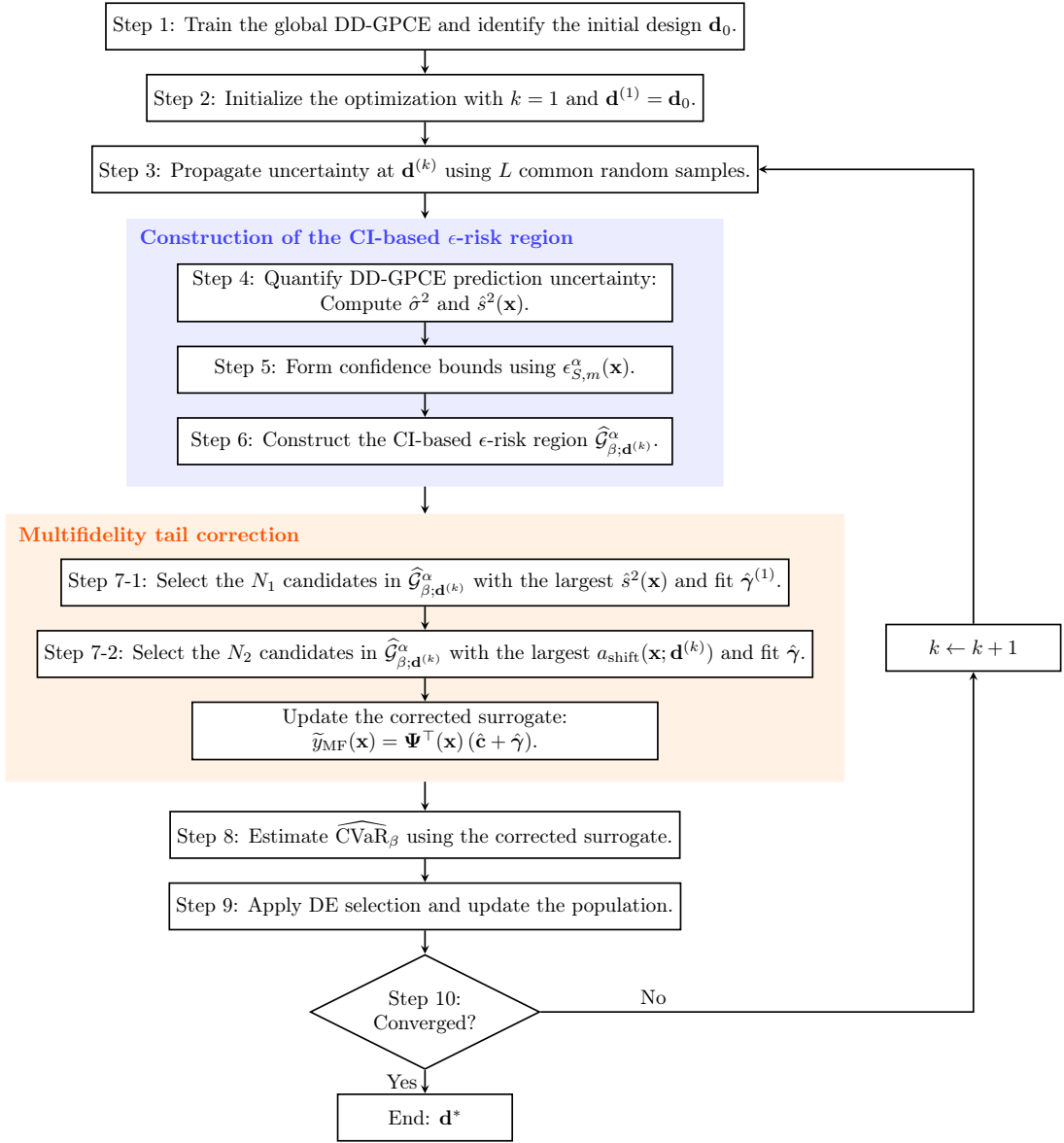
\begin{figure}[htb!]
\centering
\resizebox{0.94\textwidth}{!}{
	\begin{tikzpicture}[
		node distance=0.4cm and 1cm,
		process/.style={
			rectangle, draw=black, thick,
			minimum width=8.5cm, minimum height=0.8cm,
			align=center, fill=white
		},
		decision/.style={
			diamond, draw=black, thick, aspect=2,
			minimum width=3cm, align=center, fill=white
		},
		io/.style={
			rectangle, draw=black, thick,
			minimum width=3cm, minimum height=0.8cm,
			align=center, fill=white
		},
		arrow/.style={thick, ->, >=stealth},
		subproc/.style={
			rectangle, draw=black, thick,
			minimum width=8cm, minimum height=0.8cm,
			align=center, fill=white
		},
		sideproc/.style={
			rectangle, draw=black, thick,
			minimum width=3cm, minimum height=0.8cm,
			align=center, fill=white
		}
		]
		
		\node (s1) [process]
		{Step 1: Train the global DD-GPCE and identify the initial design $\mathbf{d}_0$.};
		
		\node (s2) [process, below=of s1]
		{Step 2: Initialize the optimization with $k=1$ and $\mathbf{d}^{(1)}=\mathbf{d}_0$.};
		
		\node (s3) [process, below=of s2]
		{Step 3: Propagate uncertainty at $\mathbf{d}^{(k)}$ using $L$ common random samples.};
		
		\node (s4) [process, below=1.2cm of s3]
		{Step 4: Quantify DD-GPCE prediction uncertainty:\\
			Compute $\hat{\sigma}^2$ and $\hat{s}^2(\mathbf{x})$.};
		
		\node (s5) [process, below=0.4cm of s4]
		{Step 5: Form confidence bounds using $\epsilon_{S,m}^{\alpha}(\mathbf{x})$.};
		
		\node (s6) [process, below=0.4cm of s5]
		{Step 6: Construct the CI-based $\epsilon$-risk region
			$\widehat{\mathcal{G}}_{\beta;\mathbf{d}^{(k)}}^\alpha$.};
		
		\node (s7_1) [subproc, below=1.6cm of s6]
		{Step 7-1: Select the $N_1$ candidates in
			$\widehat{\mathcal{G}}_{\beta;\mathbf{d}^{(k)}}^\alpha$
			with the largest $\hat{s}^2(\mathbf{x})$ and fit
			$\hat{\boldsymbol{\gamma}}^{(1)}$.};
		
		\node (s7_2) [subproc, below=0.4cm of s7_1]
		{Step 7-2: Select the $N_2$ candidates in
			$\widehat{\mathcal{G}}_{\beta;\mathbf{d}^{(k)}}^\alpha$
			with the largest $a_{\mathrm{shift}}(\mathbf{x};\mathbf{d}^{(k)})$
			and fit $\hat{\boldsymbol{\gamma}}$.};
		
		\node (s7_3) [subproc, below=0.4cm of s7_2]
		{Update the corrected surrogate:\\
			$\widetilde{y}_{\mathrm{MF}}(\mathbf{x})
			=\mathbf{\Psi}^{\top}(\mathbf{x})
			\left(\hat{\mathbf{c}}+\hat{\boldsymbol{\gamma}}\right)$.};
		
		\begin{scope}[on background layer]
			\coordinate (s4_top) at ([yshift=0.4cm]s4.north);
			\coordinate (s4_left) at ([xshift=-0.5cm]s4.west);
			\coordinate (s4_right) at ([xshift=0.5cm]s4.east);
			
			\node[
			draw=blue!80, thick,
			fill=blue!15, fill opacity=0.5,
			draw opacity=0.1,
			fit=(s4_top) (s4_left) (s4) (s4_right) (s5) (s6),
			inner sep=10pt
			] (box4_6) {};
			
			\coordinate (s7_top) at ([yshift=0.4cm]s7_1.north);
			\coordinate (s7_left) at ([xshift=-0.58cm]s7_1.west);
			\coordinate (s7_right) at ([xshift=0.58cm]s7_1.east);
			
			\node[
			draw=orange!100, thick,
			fill=orange!20, fill opacity=0.5,
			draw opacity=0.1,
			fit=(s7_top) (s7_left) (s7_1) (s7_right) (s7_2) (s7_3),
			inner sep=10pt
			] (box7) {};
		\end{scope}
		
		\node at (box4_6.north west) [
		anchor=north west,
		font=\normalsize\bfseries,
		text=blue!75,
		xshift=3pt,
		yshift=-2pt
		] {Construction of the CI-based $\epsilon$-risk region};
		
		\node at (box7.north west) [
		anchor=north west,
		font=\normalsize\bfseries,
		text=orange!65!red,
		xshift=3pt,
		yshift=-3pt
		] {Multifidelity tail correction};
		
		\node (s8) [process, below=0.5cm of box7]
		{Step 8: Estimate $\widehat{\mathrm{CVaR}}_\beta$ using the corrected surrogate.};
		
		\node (s9) [process, below=of s8]
		{Step 9: Apply DE selection and update the population.};
		
		\node (s10) [decision, below=of s9]
		{Step 10:\\Converged?};
		
		\node (end) [io, below=of s10]
		{End: $\mathbf{d}^{*}$};
		
		\node (update) [sideproc, right=0.7cm of box7]
		{$k\leftarrow k+1$};
		
		\draw [arrow] (s1) -- (s2);
		\draw [arrow] (s2) -- (s3);
		\draw [arrow] (s3) -- (box4_6.north);
		
		\draw [arrow] (s4) -- (s5);
		\draw [arrow] (s5) -- (s6);
		
		\draw [arrow] (box4_6.south) -- (box7.north);
		
		\draw [arrow] (s7_1) -- (s7_2);
		\draw [arrow] (s7_2) -- (s7_3);
		
		\draw [arrow] (box7.south) -- (s8);
		\draw [arrow] (s8) -- (s9);
		\draw [arrow] (s9) -- (s10);
		\draw [arrow] (s10) -- node[anchor=east] {Yes} (end);
		
		\draw [arrow]
		(s10.east) -| node[anchor=south, xshift=-5.5cm] {No} (update.south);
		\draw [arrow] (update.north) |- (s3.east);
		
	\end{tikzpicture}
}

\vspace{10pt}
\caption{Flowchart of the proposed tail-corrected DD-GPCE framework with DE-based risk-averse design optimization.}
\label{fig:overall_flowchart}
\end{figure}

Figure~\ref{fig:overall_flowchart} summarizes the complete procedure. The following condensed description preserves the operations needed to implement each step:

\begin{enumerate}[label=\textbf{Step \arabic*:}, leftmargin=*]
	\item \textbf{Train the global surrogate.} Sample $\mathcal{D}$, evaluate the HF model, fit the baseline DD-GPCE coefficients $\hat{\mathbf{c}}$, and identify $\mathbf{d}_0$.
	\item \textbf{Initialize optimization.} Set $k=1$ and $\mathbf{d}^{(1)}=\mathbf{d}_0$.
	\item \textbf{Propagate uncertainty.} Transform the fixed $L$ common random realizations under $\mathrm{P}_{\mathbf{d}^{(k)}}$ and evaluate the response distribution.
	\item \textbf{Quantify prediction uncertainty.} Compute $\mathbf{A}^{\top}\mathbf{A}$, $\hat{\sigma}^{2}$, and $\hat{s}^{2}(\mathbf{x})$ as defined in Section~\ref{subsec:riskregion}.
	\item \textbf{Form confidence bounds.} Evaluate the baseline response and $\epsilon_{S,m}^{\alpha}(\mathbf{x})$ to obtain lower and upper bounds for every candidate.
	\item \textbf{Construct the risk region.} Sort the lower bounds to estimate $\widehat{\mathrm{VaR}}_{\beta}$ and retain candidates whose upper bounds exceed it, following Algorithm~\ref{alg:risk_region}.
	\item \textbf{Correct the tail.} First, evaluate the $N_1$ candidates with the largest $\hat{s}^{2}$ and fit $\hat{\boldsymbol{\gamma}}^{(1)}$. Then, evaluate the $N_2$ candidates with the largest $a_{\mathrm{shift}}$ and fit the final $\hat{\boldsymbol{\gamma}}$. Combine the coefficients as $\hat{\mathbf{c}}+\hat{\boldsymbol{\gamma}}$.
	\item \textbf{Estimate CVaR.} Apply Algorithm~\ref{alg:cvar_estimation} and Eq.~\eqref{eq:cvar_estimate} to the unified corrected surrogate.
	\item \textbf{Update the DE population.} After evaluating all trial vectors in generation $k$, apply DE selection using their objective and CVaR constraint values to form the population for the next generation.
	\item \textbf{Check convergence.} Return the best feasible design if the prescribed stopping criteria are satisfied; otherwise, set $k\leftarrow k+1$ and repeat from Step~3.
\end{enumerate}

\section{Numerical examples} \label{sec:examples}

Three examples assess accuracy, scalability, and engineering applicability: the nonlinear Griewank function, a ten-bar truss, and a three-dimensional suction valve with contact and large deformation.
Examples~1 and~2 compare the proposed method with standard DD-GPCE and MF importance sampling. Standard DD-GPCE has no HF tail correction. For MF importance sampling, we follow Lee and Kramer~\cite{LeeKramer2023} but replace DD-GPCE-Kriging with DD-GPCE so that all surrogate-based methods share the same global architecture. We retain the original biasing distribution, HF allocation, and likelihood-ratio estimation procedures.

Each method evaluates CVaR repeatedly within the same DE configuration. Example~3 uses only the proposed method to demonstrate applicability to an expensive, high-dimensional finite-element problem.

For Examples~1 and~2, crude MCS with 100,000 HF samples provides reference CVaR values, while each surrogate uses 10,000 samples. A convergence study found less than $0.1\%$ change between successive estimates at this size. We compute the DD-GPCE monomial moments by QMCS using $5\times10^6$ Sobol points; these are integration points, not HF evaluations.

Python implements DE and all calculations for Examples~1 and~2. For Example~3, Python controls Abaqus/CAE 2025~\cite{AbaqusCAE2025Guide} for analysis and data extraction. All runs use a Windows workstation with two $2.80$-GHz Intel Xeon Gold 6526Y processors and 192~GB of RAM.

\subsection{Example 1: Griewank function}
\label{subsec:ex1_griewank}

\begin{figure}
    \centering
    \includegraphics[width=0.8\linewidth]{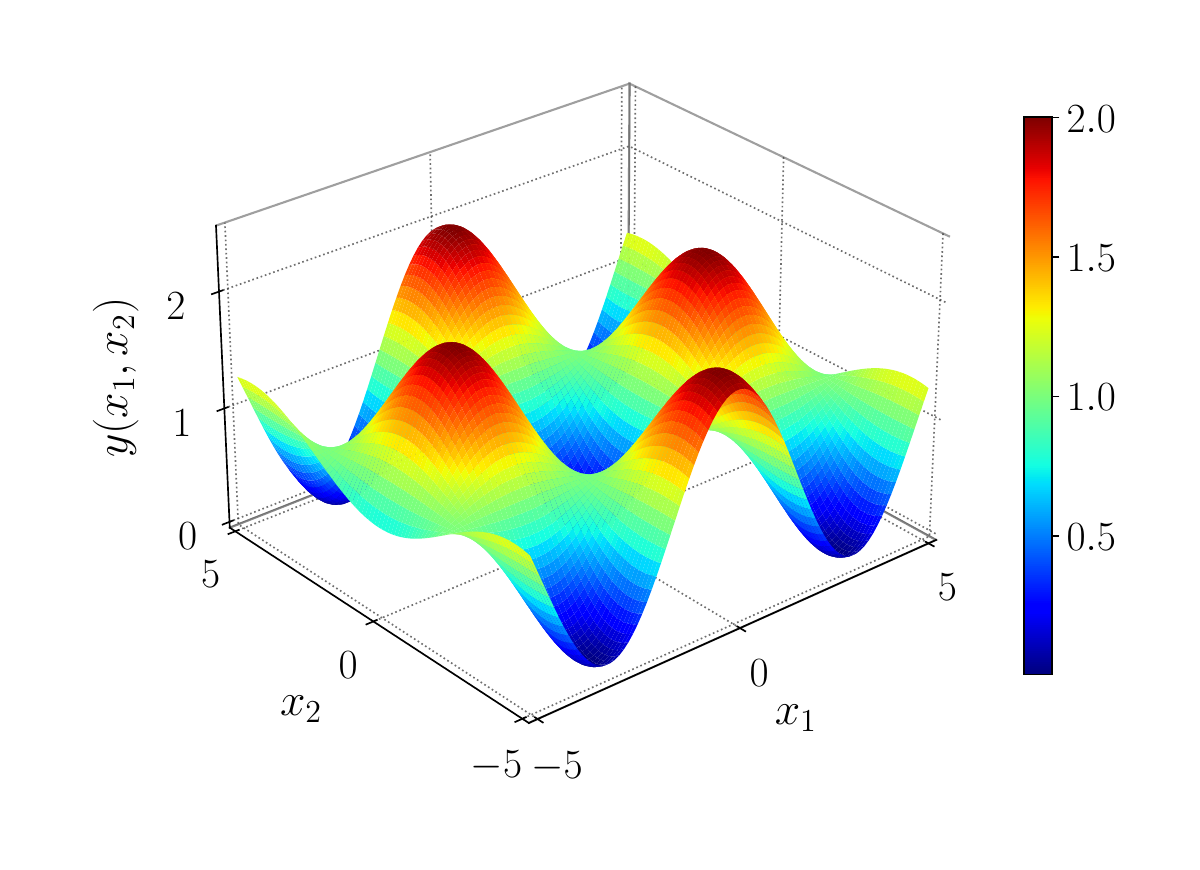}
    \caption{A graph of two-dimensional functions evaluated for $\mathrm{CVaR}_\beta$ in Example 1. The color map indicates the values of $y(x_1, x_2)$ on the domain $[-5, 5] \times [-5, 5]$.}
    \label{fig:Griewank}
\end{figure}

This example examines CVaR estimation and risk-averse design for the nonlinear, multimodal Griewank function~\cite{griewank1981generalized}. We compare the proposed method with MF importance sampling and standard DD-GPCE, using crude MCS as the reference.

\subsubsection{Problem definition}

For a bivariate random vector $\mathbf{X} = (X_1, X_2)^\top$, the true HF response is defined by the 2D Griewank function on the domain $[-5, 5]^2$:
\begin{equation}
    y(\mathbf{X}) = 1 + \frac{1}{4000} \sum_{i=1}^{2} X_i^2 - \cos(X_1)\cos\left(\frac{X_2}{\sqrt{2}}\right).
    \label{eq:griewank}
\end{equation}
Figure~\ref{fig:Griewank} shows the Griewank function on the domain $[-5, 5] \times [-5, 5]$. The function contains multiple local minima distributed throughout the domain, making it a useful benchmark for evaluating the proposed framework. We construct the multifidelity data set using a small number of exact function evaluations as HF data and surrogate-model predictions as LF data.

For the initial CVaR estimation, we model the random input vector $\mathbf{X}$ using a bivariate Gaussian distribution truncated at three standard deviations from each marginal mean. We set the mean vector to $\mathbb{E}[\mathbf{X}]=(2.0,1.0)^\top$ and the marginal standard deviations to $\sigma_1=\sigma_2=0.4$. We impose a Pearson correlation coefficient of $\rho=0.4$ between $X_1$ and $X_2$.

In the subsequent risk-averse design phase, we define the design vector as $\mathbf{d}=(d_1,d_2)^\top$, where $d_1=\mathbb{E}[X_1]$ and $d_2=\mathbb{E}[X_2]$. We optimize $\mathbf{d}$ over $\mathcal{D}=[-5,5]^2$ and initialize the search at $\mathbf{d}^{(0)}=(2.0,1.0)^\top$, corresponding to the mean vector used in the initial CVaR estimation. At each candidate design, we model $\mathbf{X}$ using a bivariate Gaussian distribution with mean $\mathbf{d}$, fixed marginal standard deviations $\sigma_1=\sigma_2=0.4$, and a Pearson correlation coefficient of $\rho=0.4$. We truncate each marginal distribution at three standard deviations from its mean.

\subsubsection{CVaR estimation results}
\label{subsubsec:ex1_cvar_results}

To evaluate the accuracy of the CVaR estimates, we compare the proposed MF tail-correction method with standard DD-GPCE and MF importance sampling, using direct MCS of the HF model as the reference. For this comparison, the same HF sample budget is assigned to MF importance sampling and the proposed method. For the proposed method, the total HF budget $L_{hf}$ is allocated between the two stages of the adaptive sampling procedure, and Tikhonov regularization with $\lambda=1.0\times10^{-6}$ is applied to the residual correction.

We use two error measures to quantify the accuracy and consistency of the CVaR estimates relative to the crude-MCS reference: the mean relative difference (MRD) and the normalized root-mean-square difference (N-RMSD). The MRD measures the average absolute relative error over $K$ independent runs:
\begin{equation}
    \text{MRD} = \frac{1}{K} \sum_{k=1}^{K} \left| \frac{\widehat{\mathrm{CVaR}}_{\beta;\mathbf{d}}^{(k)} [\hat{y}(\mathbf{X})] - \widehat{\mathrm{CVaR}}_{\beta;\mathbf{d}}[y(\mathbf{X})]}{\widehat{\mathrm{CVaR}}_{\beta;\mathbf{d}} [y(\mathbf{X})]} \right|.
    \label{eq:mrd}
\end{equation}
Here, $\widehat{\mathrm{CVaR}}_{\beta;\mathbf{d}}^{(k)}[\hat{y}(\mathbf{X})]$ denotes the surrogate-based CVaR estimate from the $k$th independent run, and $\widehat{\mathrm{CVaR}}_{\beta;\mathbf{d}}[y(\mathbf{X})]$ denotes the crude-MCS reference estimate. The N-RMSD measures the overall magnitude of the estimation errors relative to the reference:
\begin{equation}
   \text{N-RMSD} = \frac{1}{\left|\widehat{\mathrm{CVaR}}_{\beta;\mathbf{d}} [y(\mathbf{X})]\right|} \sqrt{ \frac{1}{K} \sum_{k=1}^{K} \left( \widehat{\mathrm{CVaR}}^{(k)}_{\beta;\mathbf{d}} [\hat{y}(\mathbf{X})] - \widehat{\mathrm{CVaR}}_{\beta;\mathbf{d}} [y(\mathbf{X})] \right)^2 }.
    \label{eq:nrmsd}
\end{equation}
We use $K=30$ independent runs for all evaluations.

Table~\ref{tab:cvar_griewank_results} reports results for $n_{\mathrm{hf}}=n_c=4$. Crude MCS gives the reference 1.7072, whereas uncorrected DD-GPCE gives 1.6271 and 4.56\% MRD. With four tail-focused HF evaluations, MF tail correction gives 1.7190, 2.39\% MRD, and 2.76\% N-RMSD; MF importance sampling gives 4.02\% N-RMSD. The proposed method reduces MRD by 47.6\% relative to DD-GPCE and 31.9\% relative to MF importance sampling.

\begin{table*}[htbp]
\centering
\small
\begin{threeparttable}
\caption{$\mathrm{CVaR}_{\beta}$ estimates for the Griewank function at $\beta=0.95$, obtained using the proposed MF tail-correction method, DD-GPCE-based MF importance sampling, standard DD-GPCE, and crude MCS. All surrogate-based methods use the same baseline DD-GPCE settings, $S=1$ and $m=4$.}
\label{tab:cvar_griewank_results}
\begin{tabular*}{\textwidth}{@{\extracolsep{\fill}}l c c c c c c@{}}
\toprule
Methods & $\text{CVaR}_{\beta;\mathbf{d}}$ estimate\tnote{a} & MRD\tnote{a} (\%) & N-RMSD\tnote{a} (\%) & \multicolumn{3}{c}{Model evaluations} \\
\cmidrule(lr){5-7}
 & & & & HF\tnote{b} & HF\tnote{c} & LF\tnote{d} \\
\midrule
\addlinespace
MF tail correction (Proposed) & 1.7190 & 2.3930 & 2.7606 & 100 & 4\tnote{g} & 10,000\tnote{e} \\
MF importance sampling & 1.6729 & 3.5118 & 4.0197 & 100 & 4\tnote{f} & 10,000\tnote{e} \\
standard DD-GPCE & 1.6271 & 4.5637 & 4.5646 & 100 & -- & 10,000\tnote{e} \\
crude MCS (Benchmark) & 1.7072 & -- & -- & -- & 100,000 & -- \\
\bottomrule
\end{tabular*}
\begin{tablenotes}
\scriptsize
\item[a] The estimates are averaged over $K = 30$ trials.
\item[b] The exact high-fidelity samples used to train the initial standard DD-GPCE surrogate model.
\item[c] The additional high-fidelity samples used for tail correction, importance sampling, or direct exact evaluation.
\item[d] The predictive low-fidelity output generated by the surrogate models.
\item[e] The low-fidelity output samples are used to estimate the baseline risk region and initial CVaR prior to correction.
\item[f] The exact high-fidelity output samples are injected via the MF importance sampling scheme.
\item[g] The exact high-fidelity output samples are injected for the tail correction.
\end{tablenotes}
\end{threeparttable}
\end{table*}

\begin{figure}[htbp]
    \centering
    \begin{subfigure}[b]{0.48\textwidth}
        \centering
        \includegraphics[width=\textwidth]{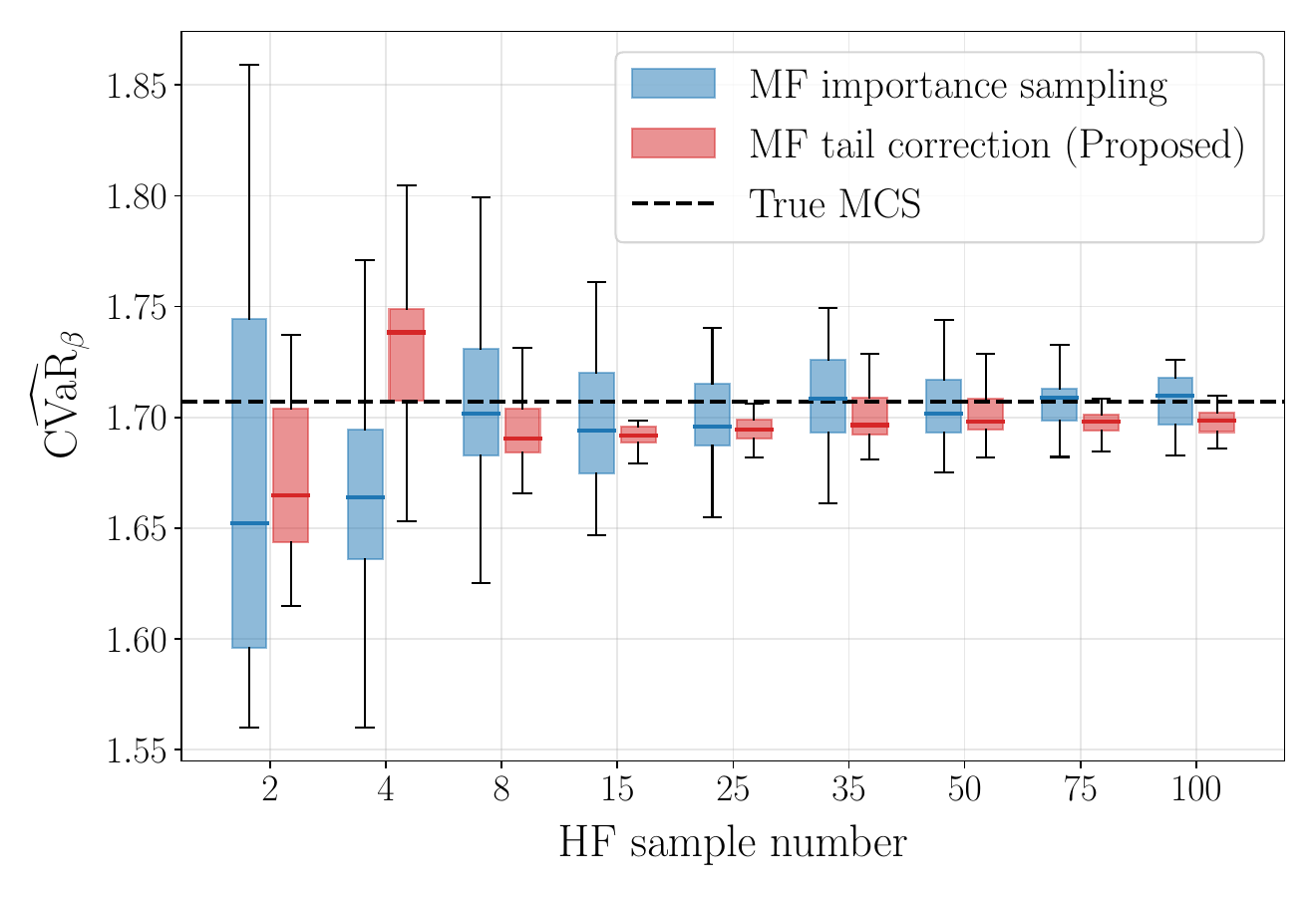}
        \caption{Distribution of CVaR estimates}
        \label{fig:cvar_trend}
    \end{subfigure}
    \hfill
    \begin{subfigure}[b]{0.48\textwidth}
        \centering
        \includegraphics[width=\textwidth]{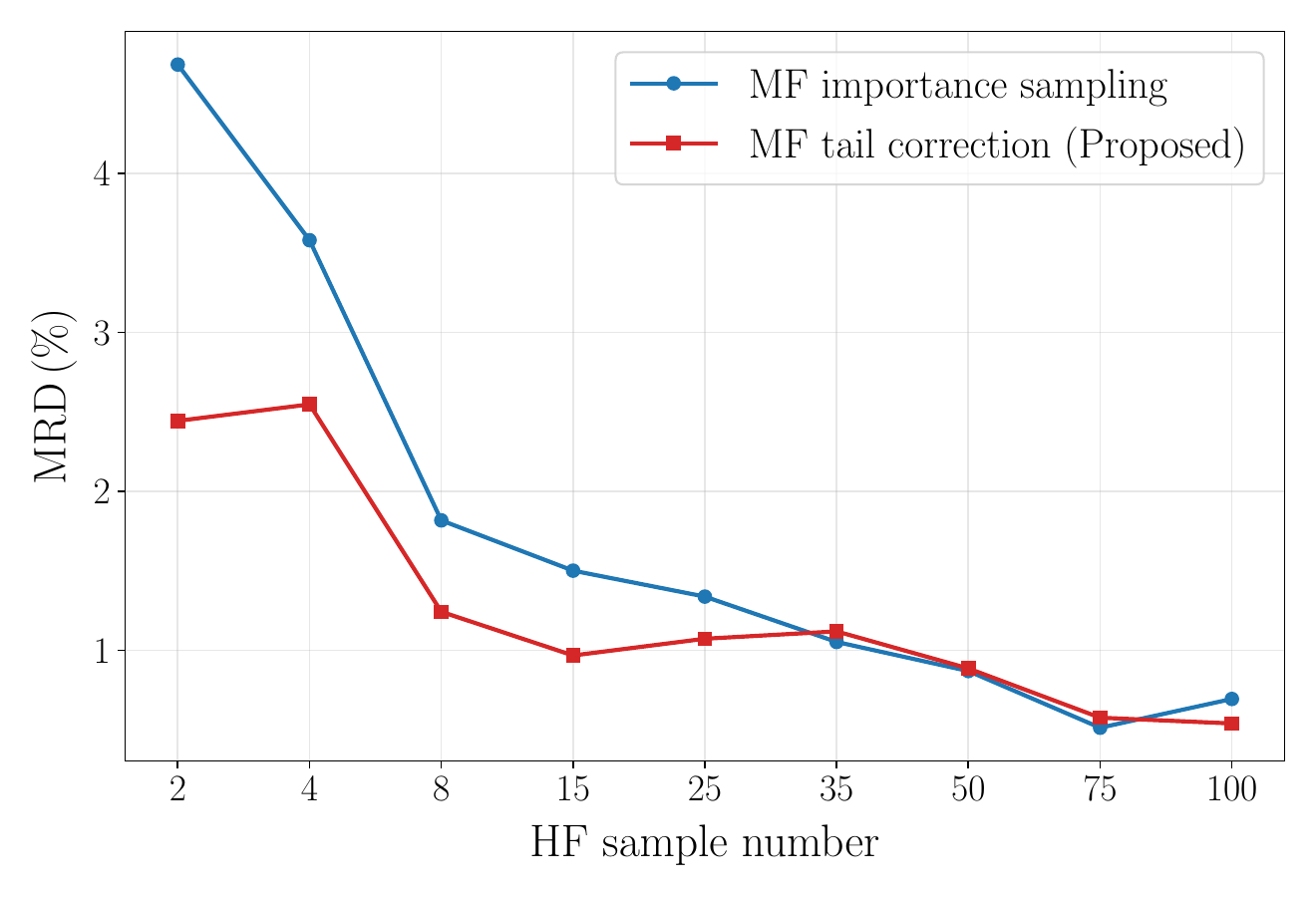}
        \caption{Mean relative difference (MRD)}
        \label{fig:mrd_trend}
    \end{subfigure}
    \caption{Griewank CVaR accuracy versus HF sample budget: (a) distributions of CVaR estimates and (b) MRD. The dashed line in (a) denotes the MCS reference.}
    \label{fig:cvar_mrd_trend}
\end{figure}

To assess accuracy and efficiency across HF budgets, Figure~\ref{fig:cvar_trend} shows the CVaR distributions and Figure~\ref{fig:mrd_trend} shows MRD.
MF importance sampling varies more at very small budgets, whereas MF tail correction has a narrower distribution and lower MRD for $L_{hf}\leq25$. At $L_{hf}=8$, the proposed method is already below 1.5\% MRD; MF importance sampling reaches this level at $L_{hf}=25$. Beyond 35 HF samples, both approach the MCS reference. This low-budget accuracy supports fewer HF evaluations per optimization generation.

\subsubsection{Risk-averse design results}
\label{subsubsec:ex1_RADO_results}

This section presents the risk-averse design results for the Griewank function. We formulate the optimization problem to minimize the deterministic Griewank response $y(\mathbf{d})$ subject to a CVaR constraint on the uncertain response $y(\mathbf{X})$:
\begin{equation}
    \begin{aligned}
        & \min\limits_{\mathbf{d} \in \mathcal{D}} & & c_0(\mathbf{d}) := y(\mathbf{d}) \\
        & \text{subject to} & & c_1(\mathbf{d}) := \widehat{\mathrm{CVaR}}_{\beta;\mathbf{d}}[y(\mathbf{X})]-0.44 \le 0, \\
        & & & -5.0 \le d_i \le 5.0, \quad i=1, 2.
    \end{aligned}
    \label{eq:griewank_RADO}
\end{equation}
The design vector $\mathbf{d}=(d_1,d_2)^\top$ specifies $\mathbb{E}_{\mathbf{d}}[X_i]=d_i$. We set the CVaR risk level to $\beta=0.95$ and the confidence-interval parameter to $\alpha=0.05$; the remaining aleatoric uncertainty follows the initial problem definition.

To compare optimization accuracy and computational efficiency, we evaluate the CVaR constraint throughout the optimization using four approaches: standard DD-GPCE, MF importance sampling, the proposed MF tail correction, and crude MCS as the reference. For consistency, all three surrogate-based approaches use the same baseline DD-GPCE model, constructed with $N=2$, $S=2$, and $m=9$.

Table~\ref{tab:RADO_griewank} summarizes the optimization results. Standard DD-GPCE uses 150 initial HF samples across the design space; MF tail correction and MF importance sampling begin with 60. After each CVaR estimate, both adaptive methods append their acquired tail-correction or importance-sampling points to the global training set and update the baseline surrogate.

After 100 generations, standard DD-GPCE returns $\mathbf{d}^{\ast}=(-0.0173,-0.1854)^\top$. It uses only the 150 initial HF evaluations because the surrogate is not updated, but its objective value is 0.0087 and $c_1(\mathbf{d}^{\ast})=0.0094>0$, so the design is infeasible.

Both adaptive methods converge near the theoretical optimum $\mathbf{d}^{\ast}=(0,0)^\top$ and satisfy the CVaR constraint, $c_1(\mathbf{d}^{\ast})=-0.0042$. MF importance sampling requires 16,230--22,610 HF evaluations, whereas MF tail correction requires 198--748, depending on $n_c$. With $n_c=4$, the proposed method uses 388 HF evaluations and yields $\Delta\mathrm{CVaR}=-5.92\times10^{-11}$.

\begin{table*}[htbp]
\centering
\renewcommand{\arraystretch}{1.25}
\caption{Risk-averse Griewank designs. Standard DD-GPCE uses 150 initial HF samples; MF tail correction and MF importance sampling each use 60.}
\label{tab:RADO_griewank}
\resizebox{\textwidth}{!}{%
\begin{tabular}{l c c c c c c c c}
\toprule
 & \multicolumn{3}{c}{Proposed MF tail correction} & \multicolumn{3}{c}{MF importance sampling} & {Standard} & Crude \\
\cmidrule(lr){2-4} \cmidrule(lr){5-7}
Metrics & $n_c=2$ & $n_c=4$ & $n_c=8$ & $n_{\mathrm{hf}}=15$ & $n_{\mathrm{hf}}=30$ & $n_{\mathrm{hf}}=50$ & DD-GPCE & MCS \\
\midrule
$d_1^*$ & $6.38 \times 10^{-6}$ & $-1.01 \times 10^{-8}$ & $-7.74 \times 10^{-9}$ & $-3.23 \times 10^{-3}$ & $1.89 \times 10^{-5}$ & $-1.96 \times 10^{-8}$ & $-0.0173$ & $0.0000^{\mathrm{a}}$ \\
$d_2^*$ & $2.49 \times 10^{-7}$ & $-1.22 \times 10^{-8}$ & $-2.29 \times 10^{-8}$ & $-9.50 \times 10^{-4}$ & $7.36 \times 10^{-5}$ & $-1.00 \times 10^{-8}$ & $-0.1854$ & $0.0000^{\mathrm{a}}$ \\
\addlinespace
{Objective $c_0(\mathbf{d}^{\ast})$} & $0.0000^{\mathrm{a}}$ & $0.0000^{\mathrm{a}}$ & $0.0000^{\mathrm{a}}$ & $5.45 \times 10^{-6}$ & $1.54 \times 10^{-9}$ & $0.0000^{\mathrm{a}}$ & $0.0087$ & $0.0000^{\mathrm{a}}$ \\
{Constraint $c_1(\mathbf{d}^{\ast})$} & $-0.0042$ & $-0.0042$ & $-0.0042$ & $-0.0041$ & $-0.0042$ & $-0.0042$ & $0.0094$ & $-0.0042$ \\
$\text{CVaR}_\beta[y(\mathbf{X})]$ & 0.4358 & 0.4358 & 0.4358 & 0.4359 & 0.4358 & 0.4358 & 0.4494 & 0.4358 \\
$\Delta \text{CVaR}^{\mathrm{b}}$ & $4.53 \times 10^{-8}$ & $-5.92 \times 10^{-11}$ & $-3.10 \times 10^{-11}$ & $-4.57 \times 10^{-5}$ & $4.72 \times 10^{-8}$ & $-1.30 \times 10^{-10}$ & $1.36 \times 10^{-2}$ & -- \\
\addlinespace
{No. of HF evaluations} & 198 & 388 & 748 & 21,585 & 16,230 & 22,610 & 150 & 4,280,000 \\
{No. of generations} & 69 & 82 & 86 & 115 & 37 & 33 & 100 & 34 \\
\bottomrule
\end{tabular}%
}
\par\vspace{1pt}
{\raggedright\footnotesize $^{\mathrm{a}}$ Values less than $1 \times 10^{-10}$ are denoted as 0.0000.\par}
{\raggedright\footnotesize $^{\mathrm{b}}$ $\Delta\mathrm{CVaR}=\widehat{\mathrm{CVaR}}_{\beta;\mathbf{d}^{\ast}}[\hat{y}(\mathbf{X})]-\widehat{\mathrm{CVaR}}_{\beta;\mathbf{d}^{\ast}}[y(\mathbf{X})]$ is the signed estimation error; $\hat{y}$ is the surrogate prediction and $y$ is the HF response evaluated by crude MCS.\par}
\end{table*}

\subsection{Example 2: Ten-bar truss structure}
\label{subsec:ex2_truss}

\begin{figure}[htbp]
    \centering
    \includegraphics[width=0.8\linewidth]{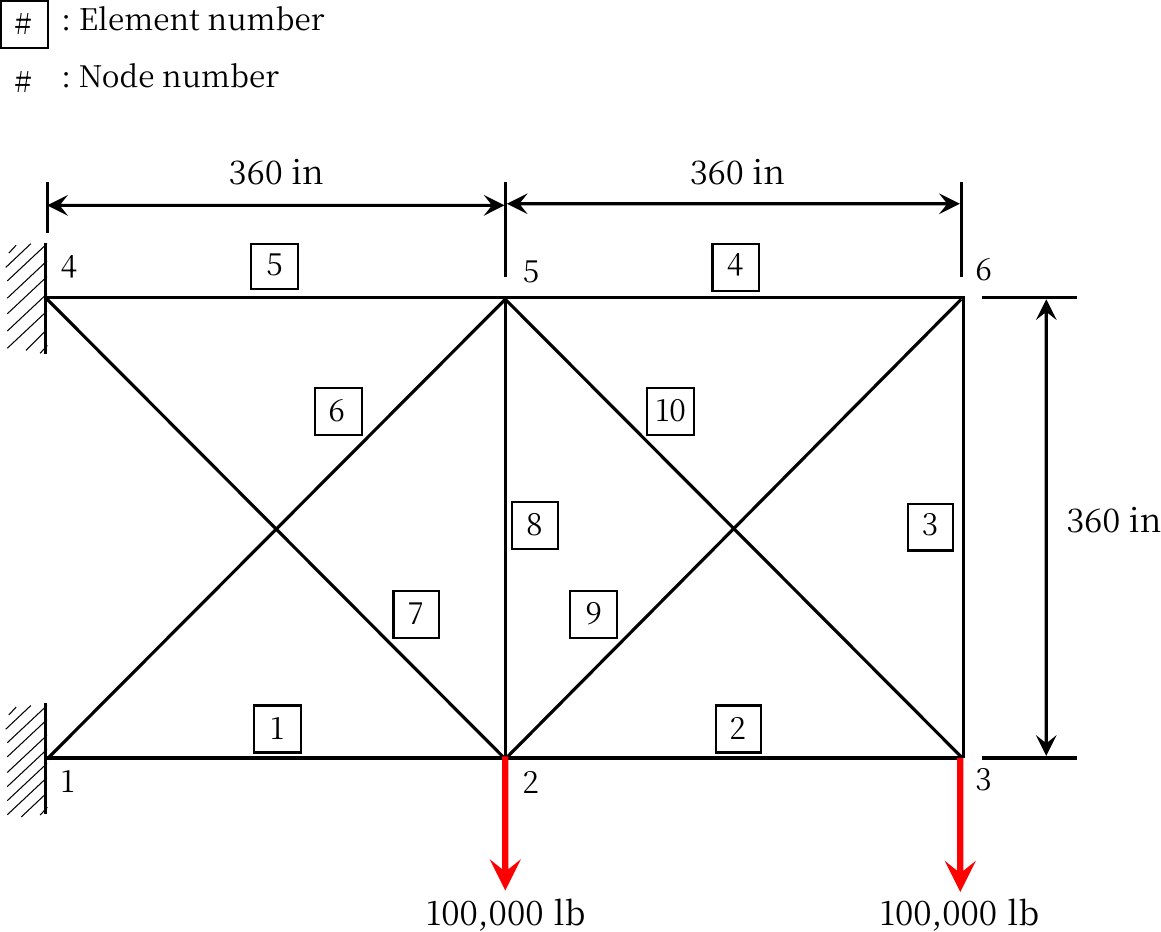}
    \caption{A ten-bar truss structure (Example 2).}
    \label{fig:truss}
\end{figure}

We next consider the ten-bar truss, a standard benchmark in structural optimization~\cite{venkayya1971design}. We use this problem to assess the proposed MF tail-correction framework in a higher-dimensional structural-mechanics setting.

Figure~\ref{fig:truss} shows that the ten-bar truss is simply supported at nodes 1 and 4, and subjected to vertically downward concentrated loads of 100,000 lb at nodes 2 and 3. The truss is composed of an aluminum alloy with Young's modulus $10^7$ psi and mass density 0.1 lb/in$^3$. The stochastic model comprises eleven ($N = 11$) random variables, $\mathbf{X} = (X_1, \dots, X_{11})^\top$. The variable $X_1$ represents the common length of the truss elements, while $X_2$ through $X_{11}$ represent the respective cross-sectional areas of the ten bars. Each variable $X_i$ follows a lognormal distribution with mean $\mathbb{E}_{\mathbf{d}}[X_i]$ and standard deviation $0.05\mathbb{E}_{\mathbf{d}}[X_i]$. Dependence among the input variables is represented by a pairwise correlation coefficient of  0.3973 for all $i\neq j$.

The design vector $\mathbf{d} = (d_1, \dots, d_{11})^\top$ is directly defined as the mean values of the eleven random variables, such that $d_i = \mathbb{E}_{\mathbf{d}}[X_i]$ for $i = 1, \dots, 11$. For the optimization process, the initial design point is set to $d_1^{(0)} = 360 \text{ in}$ and $d_i^{(0)} = 30 \text{ in}^2$ for $i = 2, \dots, 11$. The objective is to minimize the total volume of the truss while satisfying CVaR constraints on the structural response. Specifically, the upper-tail CVaR of the downward displacement at node 3, $v_3(\mathbf{X})$, must not exceed  1.1766~in, and that of the stress in bar 1, $\sigma_1(\mathbf{X})$, must not exceed 6200~psi. The resulting risk-averse design problem is 
\begin{equation}
    \begin{aligned}
        & \min\limits_{\mathbf{d} \in \mathcal{D}} & & c_0(\mathbf{d}) := \rho\int_{\mathcal{D}'(\mathbf{d})}\mathrm{d}V, \\
        & \text{subject to} & & {c_1(\mathbf{d}) := \mathrm{CVaR}_{\beta;\mathbf{d}}[v_3(\mathbf{X})]-1.1766~\mathrm{in} \le 0,} \\
        & & & {c_2(\mathbf{d}) := \mathrm{CVaR}_{\beta;\mathbf{d}}[\sigma_1(\mathbf{X})]-6200~\mathrm{psi} \le 0,} \\
        & & & 300 \text{ in} \le d_1 \le 420 \text{ in}, \\
        & & & 20 \text{ in}^2 \le d_i \le 40 \text{ in}^2, \quad i=2, \dots, 11.
    \end{aligned}
    \label{eq:tenbar_RADO}
\end{equation}
where $c_0(\mathbf{d})$ denotes the deterministic objective function for the total volume and $\beta$ is the target risk level for CVaR assessment.

\begin{figure}
    \centering
    \begin{subfigure}[b]{0.34\textwidth}
        \centering
        \includegraphics[width=\textwidth]{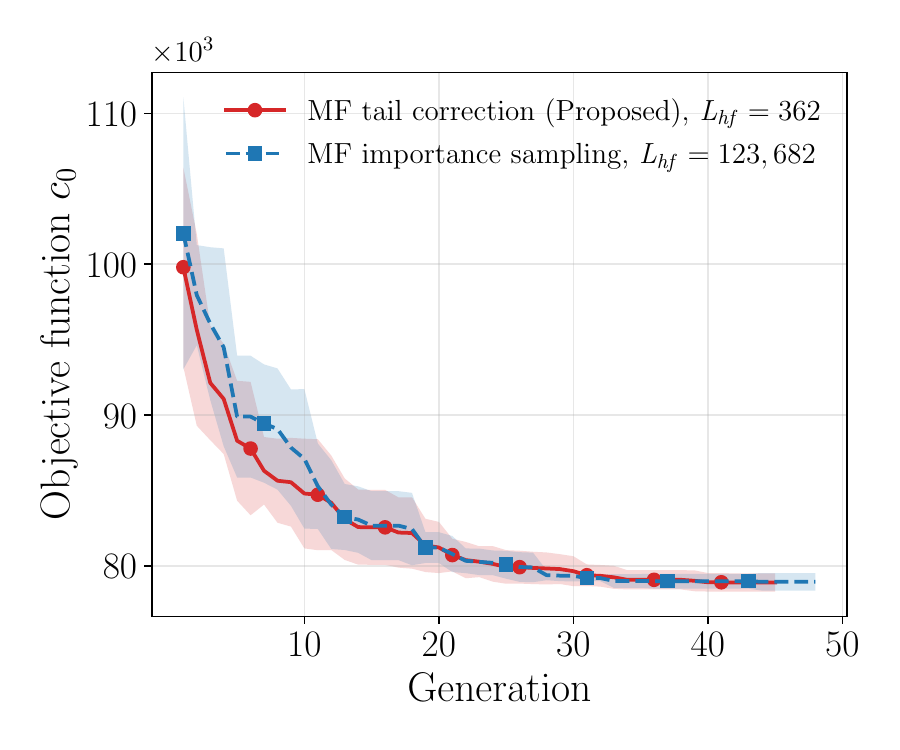}
        \caption{Objective function}
        \label{fig:conv_obj}
    \end{subfigure}
    \hfill
    \begin{subfigure}[b]{0.32\textwidth}
        \centering
        \includegraphics[width=\textwidth]{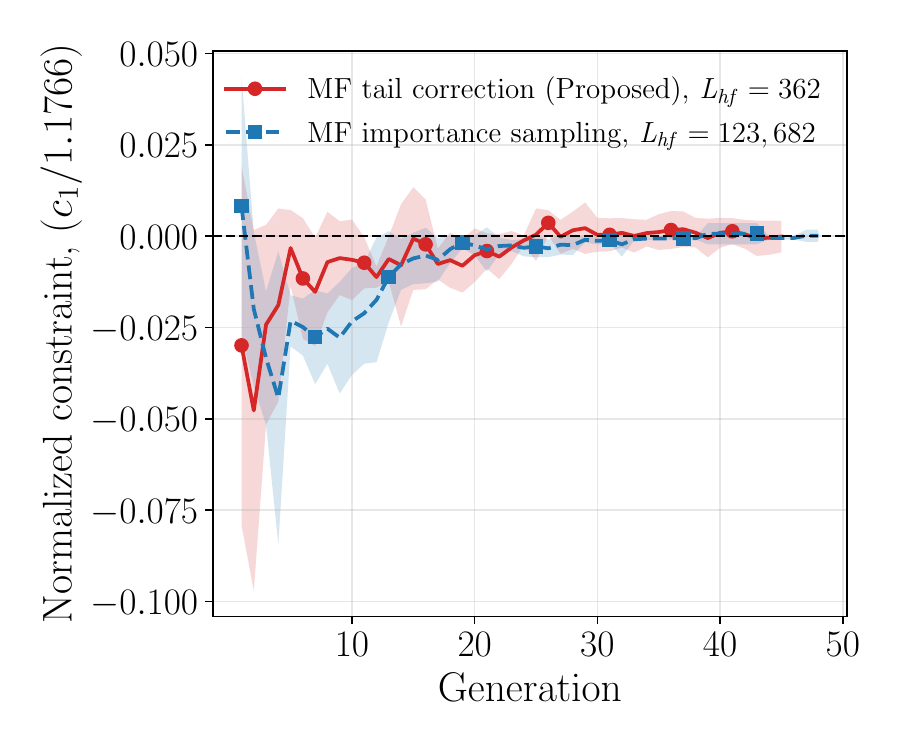}
        \caption{Normalized displacement constraint}
        \label{fig:conv_disp}
    \end{subfigure}
    \hfill
    \begin{subfigure}[b]{0.32\textwidth}
        \centering
        \includegraphics[width=\textwidth]{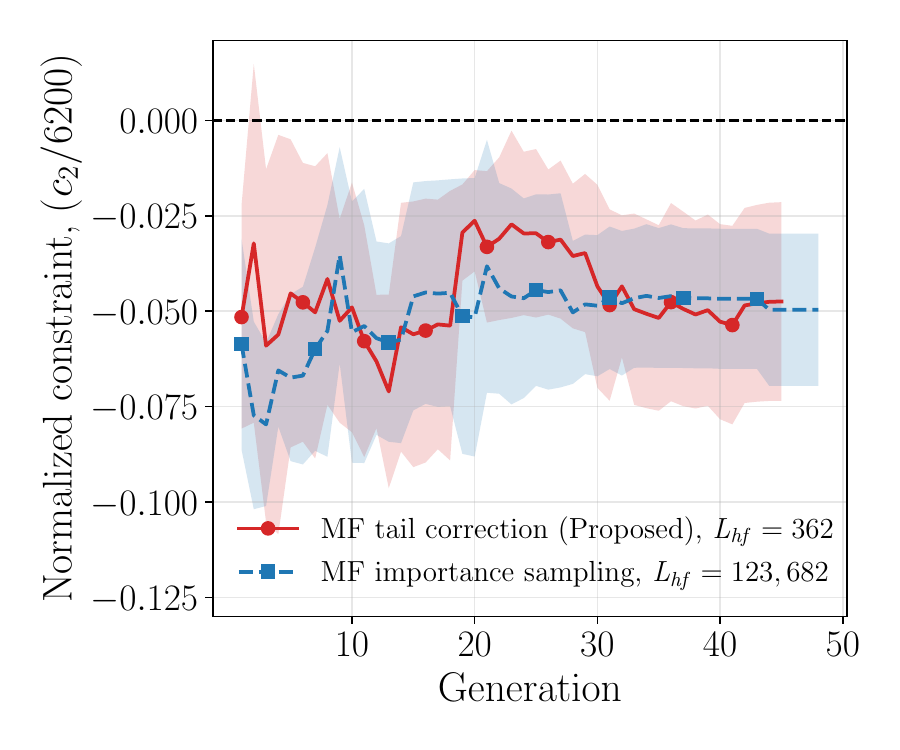}
        \caption{Normalized stress constraint}
        \label{fig:conv_stress}
    \end{subfigure}
    \caption{Convergence histories of the ten-bar truss optimization using the proposed MF tail correction and MF importance sampling methods. The solid and dashed lines represent the mean values across multiple independent trials, while the shaded regions indicate the corresponding standard deviations. (a) Objective function, (b) Normalized displacement constraint, and (c) Normalized stress constraint. The total number of high-fidelity evaluations ($L_{\mathrm{hf}}$) required by each method is reported in the legends.}
    \label{fig:tenbar_convergence}
\end{figure}

\begin{table}[htbp]
\centering
\scriptsize
\renewcommand{\arraystretch}{1.25} 
\begin{threeparttable}
\caption{{Risk-averse ten-bar truss designs obtained using MF tail correction, DD-GPCE-based MF importance sampling, standard DD-GPCE, and crude MCS.}}
\label{tab:tenbar_optimization_results}
\setlength{\tabcolsep}{4.8pt}
\begin{tabular*}{\textwidth}{@{\extracolsep{\fill}}l c c c c c c c c@{}}
\toprule
 & \multicolumn{3}{c}{Proposed MF tail correction} & \multicolumn{3}{c}{MF importance sampling} & \multicolumn{1}{c}{{Standard}} & \multicolumn{1}{c}{Crude} \\
\cmidrule(lr){2-4} \cmidrule(lr){5-7}
Metrics & $n_c=4$ & $n_c=8$ & $n_c=12$ & $n_\mathrm{hf}=15$ & $n_\mathrm{hf}=30$ & $n_\mathrm{hf}=50$ & \multicolumn{1}{c}{DD-GPCE} & \multicolumn{1}{c}{MCS} \\
\midrule
$d_1^*$, in & 300.0411 & 300.2063 & 300.2249 & 300.4593 & 300.0293 & 300.1947 & 300.0480 & 300.2050 \\
$d_2^*$, in$^2$ & 39.6062 & 38.4448 & 36.1133 & 34.1436 & 39.0487 & 37.4667 & 38.7272 & 39.3163 \\
$d_3^*$, in$^2$ & 20.6225 & 20.0069 & 20.6730 & 23.8094 & 23.3247 & 20.2680 & 20.9864 & 20.2237 \\
$d_4^*$, in$^2$ & 20.2042 & 20.8471 & 21.3286 & 20.8536 & 20.9390 & 20.0175 & 20.0962 & 21.6514 \\
$d_5^*$, in$^2$ & 20.5228 & 20.2108 & 20.0034 & 20.4187 & 20.8819 & 21.7895 & 20.1076 & 20.0858 \\
$d_6^*$, in$^2$ & 35.0545 & 39.6922 & 36.9524 & 38.4430 & 35.5132 & 36.8739 & 38.3158 & 35.7728 \\
$d_7^*$, in$^2$ & 31.7143 & 34.2496 & 39.5485 & 39.3639 & 26.7973 & 33.2775 & 29.0277 & 36.8960 \\
$d_8^*$, in$^2$ & 25.1248 & 26.2477 & 22.8250 & 21.3742 & 25.3014 & 27.6541 & 20.5646 & 26.0287 \\
$d_9^*$, in$^2$ & 20.0041 & 20.0034 & 20.4729 & 20.2258 & 20.0326 & 20.3143 & 20.5670 & 20.9761 \\
$d_{10}^*$, in$^2$ & 20.8478 & 20.5030 & 20.2471 & 20.1338 & 21.5654 & 22.3410 & 20.9622 & 22.3243 \\
$d_{11}^*$, in$^2$ & 27.9907 & 20.7254 & 26.5381 & 26.7126 & 31.7617 & 23.2503 & 26.4978 & 22.0026 \\
\addlinespace
Objective $c_0(\mathbf{d}^*)$ & 78,518.29 & 78,333.12 & 79,470.15 & 79,765.49 & 79,557.56 & 79,027.07 & 76,768.00 & 79,637.68 \\
Constraint $c_1(\mathbf{d}^*)$ & 0.0034 & 0.0025 & $-0.0014$ & 0.0017 & 0.0033 & 0.0006 & 0.0414 & $-0.0005$ \\
Constraint $c_2(\mathbf{d}^*)$ & $-441.94$ & $-348.39$ & $-284.40$ & $-81.70$ & $-224.34$ & $-94.99$ & $-510.31$ & $-488.74$ \\
$\mathrm{CVaR}_{\beta;\mathbf{d}^*}[v_3(\mathbf{X})]$ & 1.1800 & 1.1791 & 1.1752 & 1.1783 & 1.1799 & 1.1772 & 1.2180 & 1.1761 \\
$\mathrm{CVaR}_{\beta;\mathbf{d}^*}[\sigma_1(\mathbf{X})]$ & 5758.06 & 5851.61 & 5915.60 & 6118.31 & 5975.66 & 6105.01 & 5689.69 & 5711.26 \\
\addlinespace
{No. of HF evaluations} & 322 & 470 & 618 & 44,698 & 128,687 & 169,462 & 150 & 10,680,000 \\
{No. of Generations} & 43 & 40 & 39 & 25 & 37 & 30 & 39 & 18 \\
\bottomrule
\end{tabular*}
\end{threeparttable}
\end{table}

Table~\ref{tab:tenbar_optimization_results} compares the methods with a 10,680,000-sample crude-MCS reference. All MF importance sampling configurations violate the displacement constraint, with HF CVaR estimates of 1.1772--1.1799~in. MF tail correction with $n_c=12$ is feasible and gives $79{,}470.15~\mathrm{in}^3$, close to the MCS value $79{,}637.68~\mathrm{in}^3$, using 618 HF evaluations rather than as many as 169,462. Standard DD-GPCE uses only 150 initial evaluations but gives a displacement CVaR of 1.2180~in. Figure~\ref{fig:tenbar_convergence} shows that the proposed method achieves convergence trends similar to MF importance sampling with far fewer HF evaluations.

\subsection{Example 3: Suction valve system}
\label{subsec:ex3_valve}

The third example considers a suction valve system comprising a reed-type valve and a piston assembly. Reed valves are widely used in hermetic compressors as pressure-actuated components that control the flow of the working fluid, and their dynamic response directly influences the efficiency and reliability of the compressor~\cite{egger2020multi}. Mounted on the upper surface of the piston, the valve regulates fluid inflow through pressure-induced deformation.

\begin{figure}
    \centering
    \includegraphics[width=0.85\linewidth]{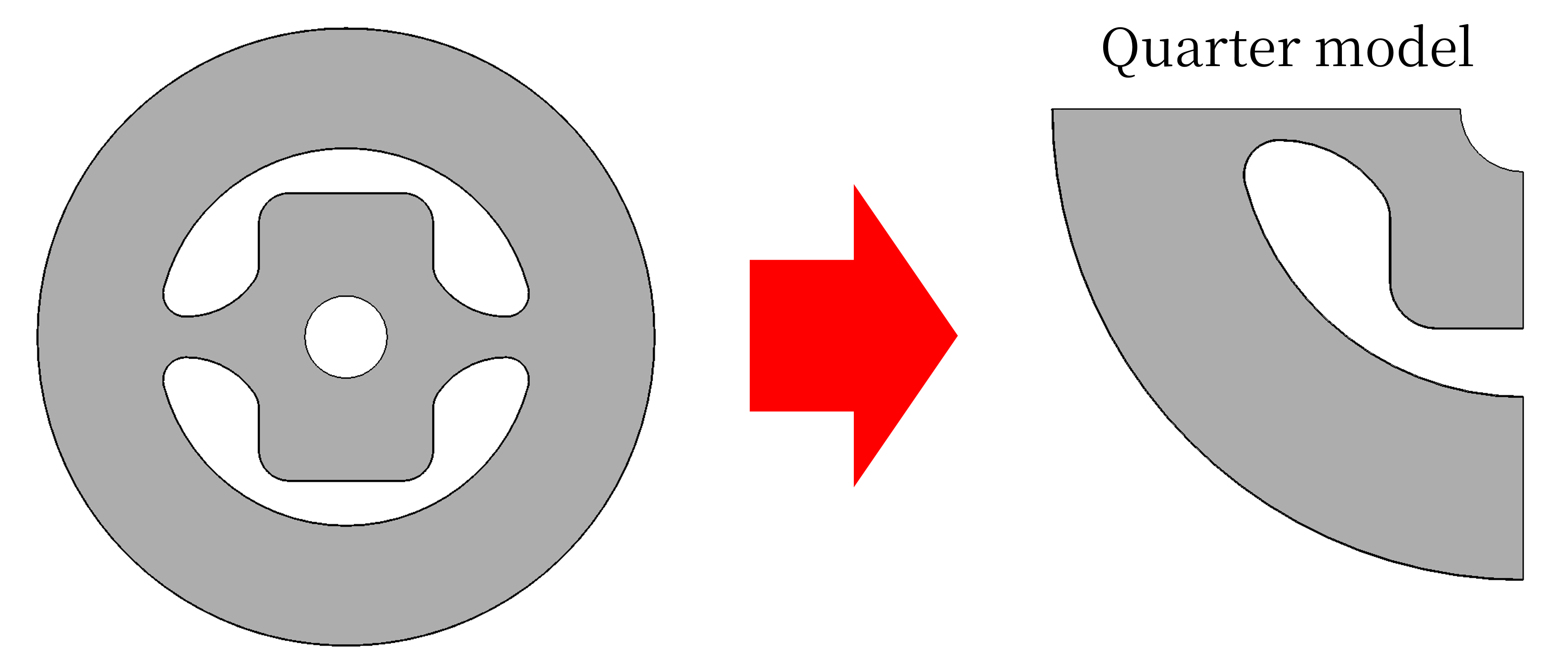}
    \caption{Valve and geometric domain reduction: the valve extracted from the full model and the quarter-symmetry model used for finite-element analysis.}
    \label{fig:suction_sys}
\end{figure}

Figure~\ref{fig:suction_sys} shows the valve and quarter-symmetry computational domain. The parameterized assembly retains the piston and its design variables. Confidentiality restrictions prevent disclosure of specific piston geometry, dimensions, and variable mappings, but the finite-element analyses and optimization include all geometric features.

Repeated HF analyses of the complete assembly would impose a substantial computational burden within the optimization loop. Accordingly, two orthogonal symmetry planes are introduced, allowing only one-quarter of the valve--piston assembly to be modeled. The resulting quarter model preserves the relevant pressure loading, contact interfaces, and boundary conditions of the complete system while substantially reducing the finite-element model size. This example evaluates the proposed tail-corrected DD-GPCE framework for structural design under geometric uncertainty, contact interactions, and pressure-induced deformation. Only the relevant non-proprietary information is reported.

\subsubsection{System description and problem definition}
The structural response of the suction valve was evaluated under two loading cases, hereafter referred to as the bending and bulging modes. In the bending mode, the suction-induced pressure difference deflects the valve away from the piston surface, thereby opening the flow passage. The resulting maximum von Mises stress is defined as the bending-stress response, $y_{\mathrm{bending}}(\mathbf{X})$. In the bulging mode, the pressure acts in the opposite direction and forces the closed valve against the piston contact surface, producing localized deformation and stress concentrations. The corresponding maximum von Mises stress is defined as the bulging-stress response, $y_{\mathrm{bulging}}(\mathbf{X})$. The two modes were evaluated using separate finite-element analyses because they involve opposing loading directions and distinct contact states. The valve opening height, denoted by $h(\mathbf{X})$, was also extracted from the bending analysis as a measure of the functional performance of the suction system.

The risk-averse design problem seeks to minimize the total mass of the system while limiting the upper-tail stress responses and maintaining sufficient  valve-opening performance. The risk level was set to $\beta=0.99$, such that each stress constraint governs the mean response within the upper $1\%$ tail of the corresponding stress distribution. The allowable CVaR thresholds for the bending and bulging stresses were set $15\%$ and $20\%$ below their respective initial-design values and are denoted by $\sigma_{\mathrm{lim}}^{(\mathrm{bending})}$ and $\sigma_{\mathrm{lim}}^{(\mathrm{bulging})}$, respectively. In addition, the expected valve opening height was required to satisfy $\mathbb{E}[h(\mathbf{X})]\geq h_{\min}$, where $h_{\min}$ denotes the prescribed minimum opening needed to maintain the suction flow passage. The numerical values of the stress thresholds and minimum opening requirement are withheld because of confidentiality restrictions. The resulting optimization problem is 
\begin{equation}
\begin{aligned}
    \min_{\mathbf{d}\in\mathcal{D}}
    \quad
    & c_{0}(\mathbf{d})
    :=
    \int_{\mathcal{D}'(\mathbf{d})}\rho(\mathbf{x})\mathrm{d}V,
    \\[4pt]
    \text{subject to}
    \quad
    & {c_{1}(\mathbf{d})
    :=
    \operatorname{CVaR}_{\beta;\mathbf{d}}
    \!\left[
        y_{\mathrm{bending}}(\mathbf{X})
    \right]
    -\sigma_{\mathrm{lim}}^{(\mathrm{bending})}
    \leq0,}
    \\
    & {c_{2}(\mathbf{d})
    :=
    \operatorname{CVaR}_{\beta;\mathbf{d}}
    \!\left[
        y_{\mathrm{bulging}}(\mathbf{X})
    \right]
    -\sigma_{\mathrm{lim}}^{(\mathrm{bulging})}
    \leq0,}
    \\[2.5pt]
    & {c_{3}(\mathbf{d})
    :=h_{\min}-\mathbb{E}\!\left[h(\mathbf{X})\right]
    \leq0,}
    \\[2pt]
    & d_{l,L}\leq d_l\leq d_{l,U},
    \qquad l=1,\ldots,18.
\end{aligned}
\label{eq:ex3_optimization}
\end{equation}

\begin{figure}[!ht]
    \centering
    \includegraphics[width=0.95\linewidth]{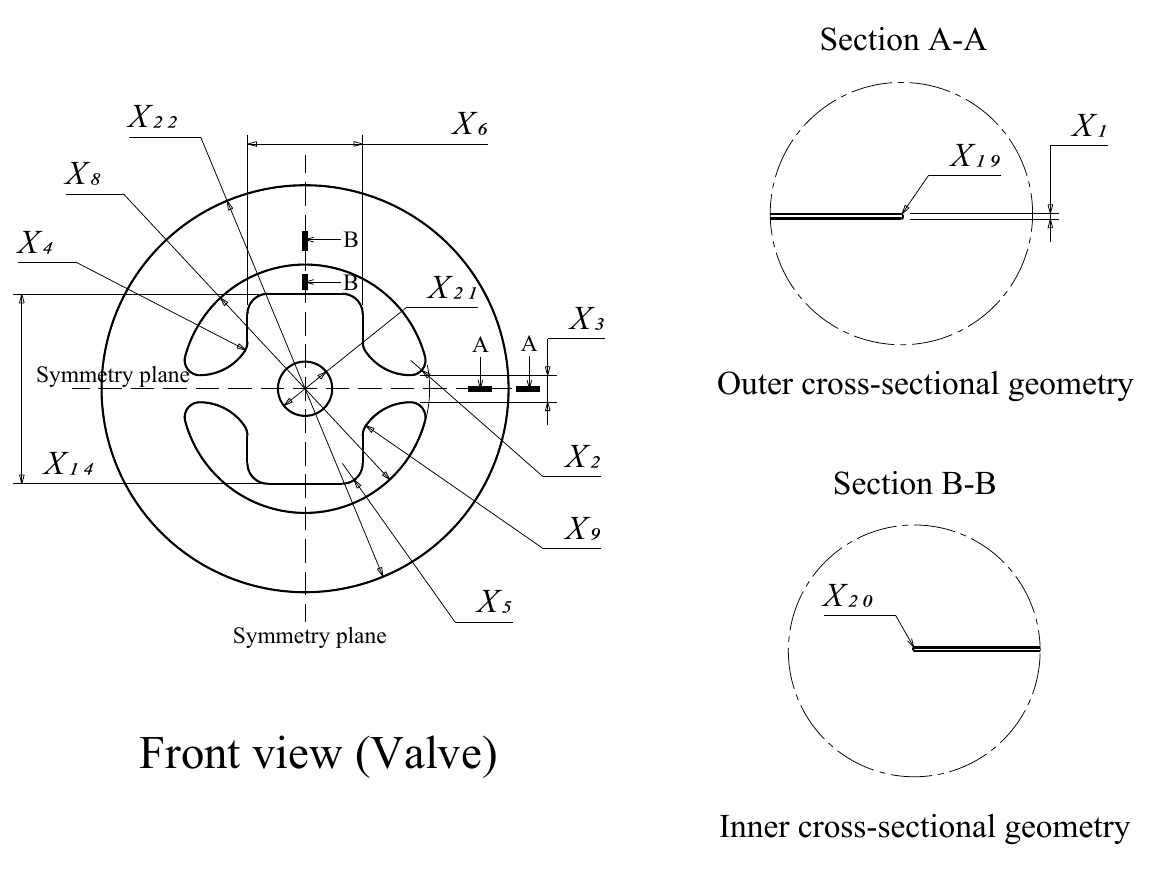}
    \caption{Parameterized valve geometry with symmetry planes and cross-sectional details}
    \label{fig:draft}
\end{figure}

The valve--piston assembly was parameterized in CATIA V5 using 24 geometric variables. Figure~\ref{fig:draft} illustrates the valve geometry and the valve-related parameters. The detailed piston geometry, dimensions, and mapping between the piston features and their associated geometric variables are omitted because of industrial confidentiality restrictions. 

We model $X_1,\ldots,X_{18}$ as random variables with design-dependent means, such that $\mathbb{E}[X_i]=d_i$ for $i=1,\ldots,18$. The variables $X_{19},\ldots,X_{24}$ capture geometric uncertainties, and we hold their mean values fixed throughout the optimization. We define the design vector as $\mathbf{d}=(d_1,\ldots,d_{18})^\top$. To protect proprietary dimensional information, we report the initial and optimal designs and their corresponding bounds in Table~\ref{tab:design_variables} as dimensionless ratios normalized by the respective initial values.

Table~\ref{tab:statistical_properties} summarizes the statistical properties of the geometric variables using relative measures.  All variables were modeled as mutually independent uniform random variables; hence, their pairwise correlation coefficients satisfy $\rho_{ij}=0$ for $i\neq j$. For $X_{1},\ldots,X_{18}$, the standard deviation was specified as $\sigma_i=c_i d_i$, where $c_i$ is the prescribed coefficient of variation. Consequently, the uncertainty interval varies with the current design value $d_i$. The variables $X_{19},\ldots,X_{24}$ were assigned fixed means and standard deviations in the computational model, but their dimensional values are not disclosed.  Instead, we characterize their uncertainty using the coefficient of variation $c_i=\sigma_i/\mu_i$.

\begin{table}[htb]
    \centering
    \begin{threeparttable}
        \caption{Normalized optimal values and prescribed bounds of the 18 design variables for the suction system}
        \label{tab:design_variables}

        \footnotesize
        \renewcommand{\arraystretch}{1.10}
        \setlength{\tabcolsep}{7pt}

        \begin{tabular}{
            @{}
            rrrr
            @{\hspace{2em}}
            rrrr
            @{}
        }
            \toprule
            $l$
            & $d_l^{*}/d_{l,0}$
            & $d_{l,L}/d_{l,0}$
            & $d_{l,U}/d_{l,0}$
            &
            $l$
            & $d_l^{*}/d_{l,0}$
            & $d_{l,L}/d_{l,0}$
            & $d_{l,U}/d_{l,0}$ \\
            \midrule

             1 & 1.0227 & 0.7727 & 1.2273
            & 10 & 1.0254 & 1.0000 & 1.0439 \\

             2 & 2.0570 & 0.7273 & 2.2727
            & 11 & 0.9640 & 0.8859 & 1.1074 \\

             3 & 0.5573 & 0.5000 & 1.2500
            & 12 & 0.9970 & 0.9591 & 1.0380 \\

             4 & 0.7853 & 0.3846 & 1.5385
            & 13 & 1.0001 & 0.9803 & 1.0185 \\

             5 & 1.5055 & 0.1333 & 1.8667
            & 14 & 0.8626 & 0.8065 & 1.0000 \\

             6 & 1.0490 & 1.0000 & 1.1765
            & 15 & 2.1981 & 0.8710 & 2.2903 \\

             7 & 0.8917 & 0.7857 & 1.2143
            & 16 & 2.1470 & 1.0000 & 3.1008 \\

             8 & 0.9982 & 0.9508 & 1.0492
            & 17 & 2.2586 & 1.0000 & 3.1008 \\

             9 & 1.1629 & 1.0000 & 1.2500
            & 18 & 1.0244 & 0.6124 & 2.9380 \\

            \bottomrule
        \end{tabular}

        \begin{tablenotes}[flushleft]
            \footnotesize
            \item All values are dimensionless and normalized by the
            corresponding initial design value $d_{l,0}$. Therefore, the
            normalized initial value of every design variable is unity.
        \end{tablenotes}
    \end{threeparttable}
\end{table}

\begin{table}[htb]
    \centering
    \begin{threeparttable}
        \caption{Statistical definitions of the 24 random variables used in the suction system}
        \label{tab:statistical_properties}

        \footnotesize
        \renewcommand{\arraystretch}{1.08}
        \setlength{\tabcolsep}{8pt}

        \begin{tabular}{@{}lllll@{}}
            \toprule
            \makecell[l]{Random\\variable}
            & Property
            & \makecell[l]{Mean\\(mm)}
            & \makecell[l]{Standard\\deviation (mm)}
            & \makecell[l]{Probability\\distribution} \\
            \midrule

            $X_{1}$  & Thickness & $d_{1}$  & $0.02d_{1}$  & Uniform \\
            $X_{2}$  & Radius    & $d_{2}$  & $0.02d_{2}$  & Uniform \\
            $X_{3}$  & Width     & $d_{3}$  & $0.02d_{3}$  & Uniform \\
            $X_{4}$  & Radius    & $d_{4}$  & $0.02d_{4}$  & Uniform \\
            $X_{5}$  & Radius    & $d_{5}$  & $0.02d_{5}$  & Uniform \\
            $X_{6}$  & Length    & $d_{6}$  & $0.02d_{6}$  & Uniform \\
            $X_{7}$  & Length    & $d_{7}$  & $0.02d_{7}$  & Uniform \\
            $X_{8}$  & Diameter  & $d_{8}$  & $0.02d_{8}$  & Uniform \\
            $X_{9}$  & Radius    & $d_{9}$  & $0.02d_{9}$  & Uniform \\
            $X_{10}$ & Diameter  & $d_{10}$ & $0.01d_{10}$ & Uniform \\
            $X_{11}$ & Length    & $d_{11}$ & $0.02d_{11}$ & Uniform \\
            $X_{12}$ & Length    & $d_{12}$ & $0.02d_{12}$ & Uniform \\
            $X_{13}$ & Radius    & $d_{13}$ & $0.01d_{13}$ & Uniform \\
            $X_{14}$ & Length    & $d_{14}$ & $0.02d_{14}$ & Uniform \\
            $X_{15}$ & Length    & $d_{15}$ & $0.02d_{15}$ & Uniform \\
            $X_{16}$ & Radius    & $d_{16}$ & $0.02d_{16}$ & Uniform \\
            $X_{17}$ & Radius    & $d_{17}$ & $0.02d_{17}$ & Uniform \\
            $X_{18}$ & Radius    & $d_{18}$ & $0.02d_{18}$ & Uniform \\
            $X_{19}$ & Radius    & $\mu_{19}$
                      & $0.02\mu_{19}$ & Uniform \\
            $X_{20}$ & Radius    & $\mu_{20}$
                      & $0.02\mu_{20}$ & Uniform \\
            $X_{21}$ & Diameter  & $\mu_{21}$
                      & $0.02\mu_{21}$ & Uniform \\
            $X_{22}$ & Diameter  & $\mu_{22}$
                      & $0.02\mu_{22}$ & Uniform \\
            $X_{23}$ & Radius    & $\mu_{23}$
                      & $0.02\mu_{23}$ & Uniform \\
            $X_{24}$ & Diameter  & $\mu_{24}$
                      & $0.02\mu_{24}$ & Uniform \\

            \bottomrule
        \end{tabular}

        \begin{tablenotes}[flushleft]
            \footnotesize
            \item The variables $X_{1},\ldots,X_{18}$ have design-dependent means. For $X_{19},\ldots,X_{24}$, $\mu_i$ denotes a fixed mean whose dimensional value is not disclosed. Their standard deviations are reported relative to the corresponding fixed means. The support of each uniform distribution is defined as            $\left[\mu_i-\sqrt{3}\sigma_i,\,            \mu_i+\sqrt{3}\sigma_i\right]$.
        \end{tablenotes}
    \end{threeparttable}
\end{table}

\begin{figure}[htbp]
    \centering
    \begin{subfigure}[b]{0.55\textwidth}
        \centering
        \includegraphics[width=\linewidth]{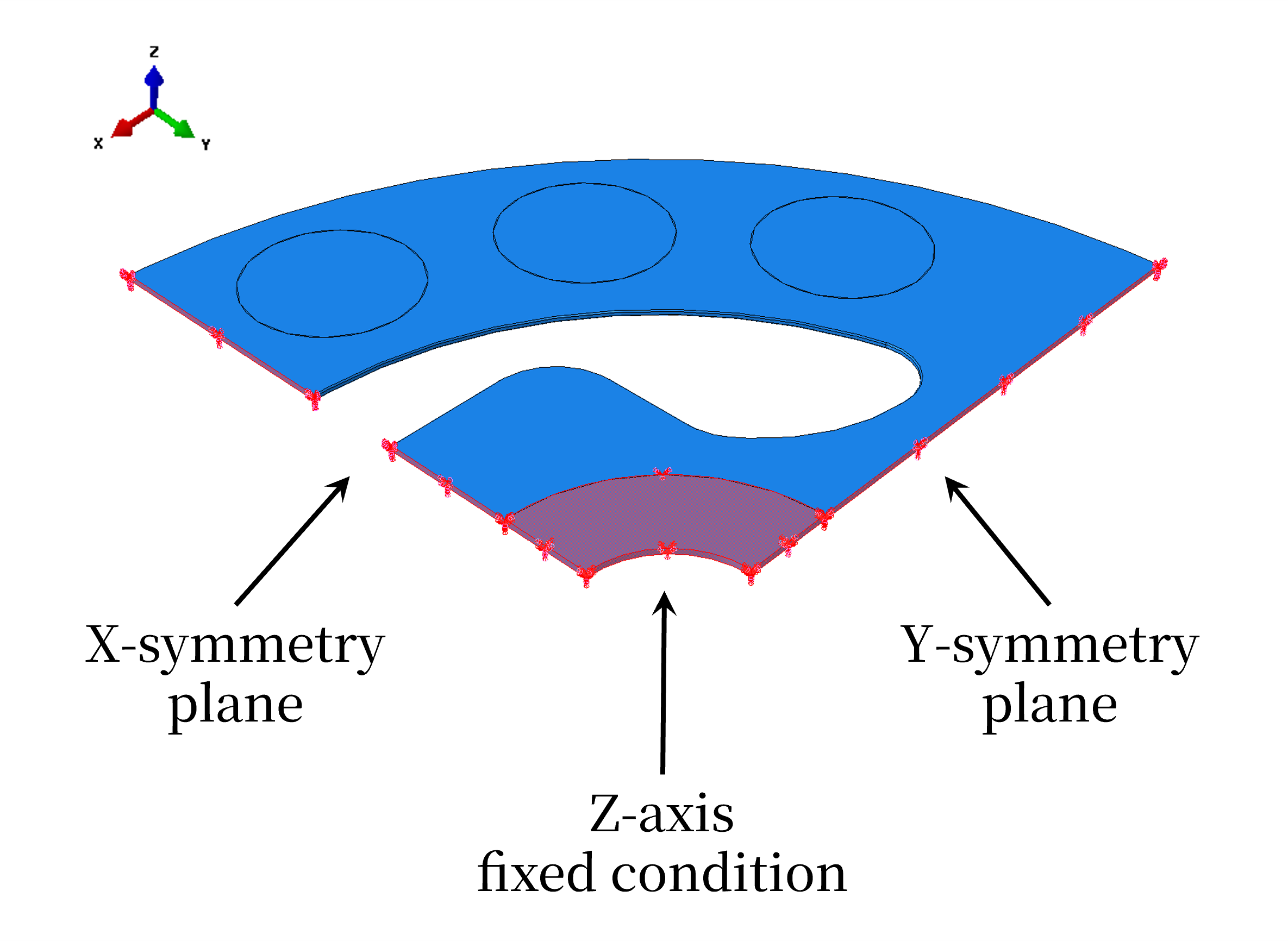}
        \caption{Symmetry and support boundary conditions of the quarter model}
        \label{fig:BC}
    \end{subfigure}
    
    \vspace{7mm}
    
    \centering
    \begin{subfigure}[b]{0.48\textwidth}
        \centering
        \includegraphics[width=\textwidth]{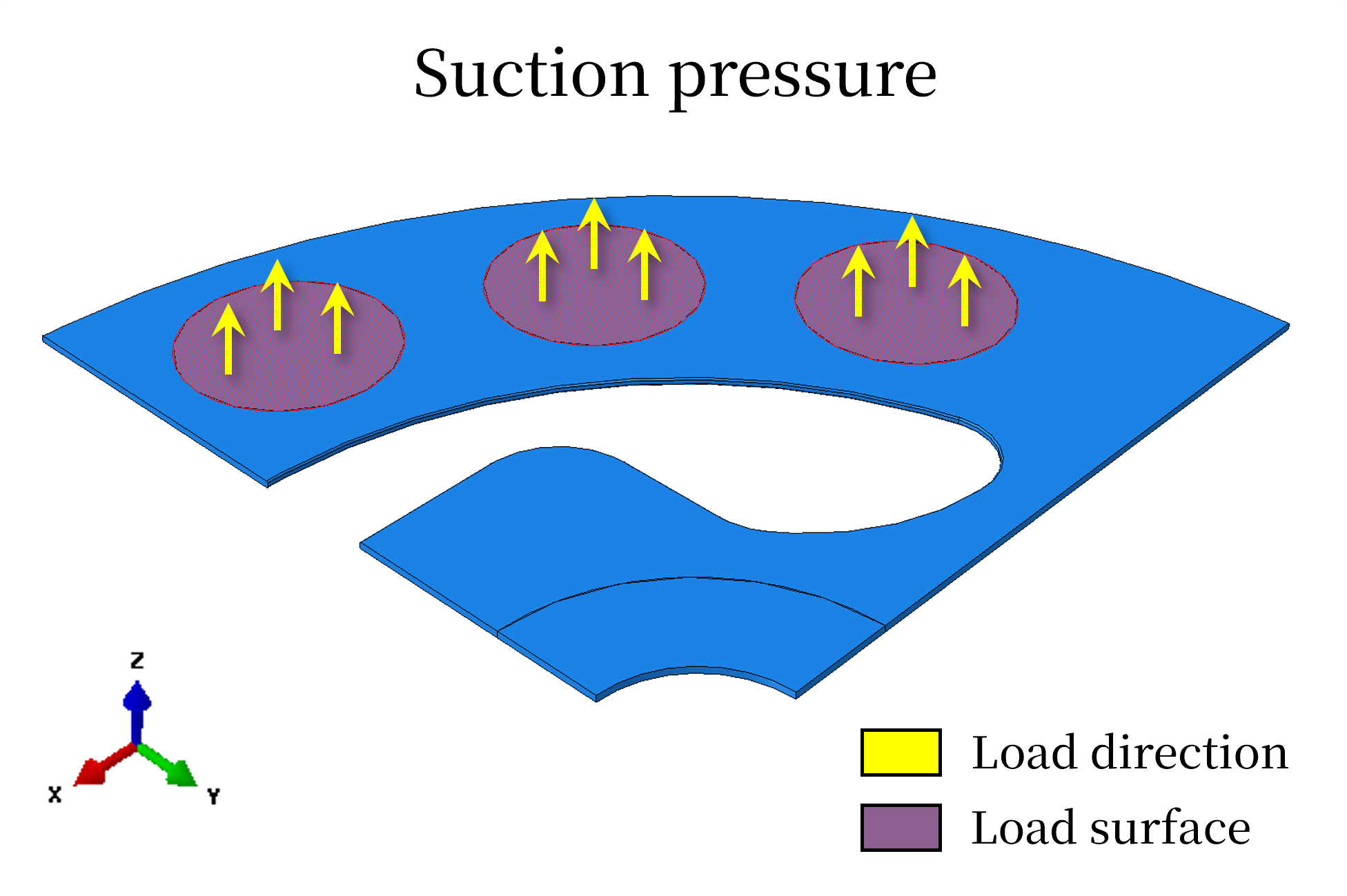}
        \caption{Prescribed pressure loading for the valve-bending condition}
        \label{fig:bending}
    \end{subfigure}
    \hfill
    \begin{subfigure}[b]{0.48\textwidth}
        \centering
        \includegraphics[width=\textwidth]{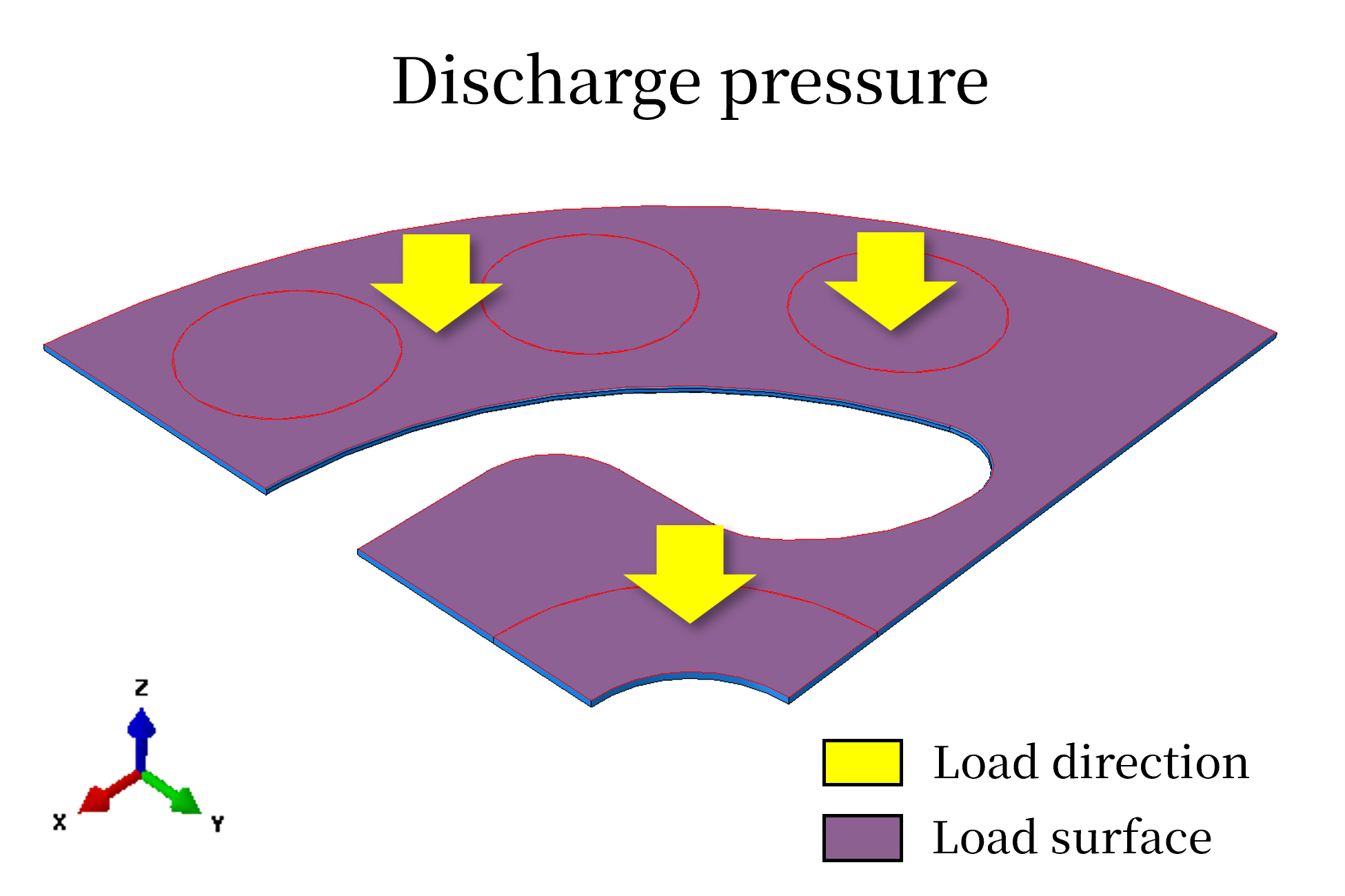}
        \caption{Prescribed pressure loading for the valve-bulging condition}
        \label{fig:bulging}
    \end{subfigure}
    \caption{Boundary conditions and prescribed pressure-loading cases for the bending and bulging analyses of the quarter model.}
    \label{fig:load}
\end{figure}

We evaluated the structural response of the quarter model using Abaqus/Explicit.  As shown in Figure~\ref{fig:BC}, we imposed symmetry boundary conditions on both cut planes by constraining the displacement component normal to each plane. We also constrained the $z$-direction displacement in the bolted region of the valve and fully fixed the piston to suppress rigid-body motion. We applied the same boundary conditions in the bending and bulging analyses.

\begin{figure}[!ht]
    \centering
    \begin{subfigure}{1.0\textwidth}
        \centering
        \includegraphics[width=\textwidth]{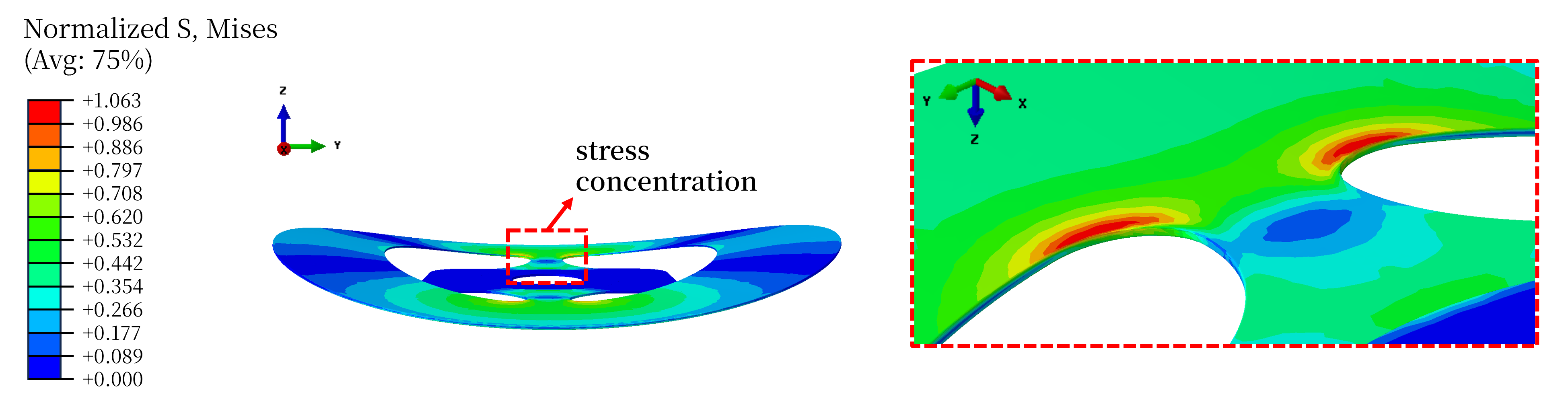}
        \caption{Normalized von Mises stress contour and deformed configuration under the bending load case.}
        \label{fig:bending_contour}
    \end{subfigure}

    \vspace{7mm}

    \begin{subfigure}{1.0\textwidth}
        \centering
        \includegraphics[width=\textwidth]{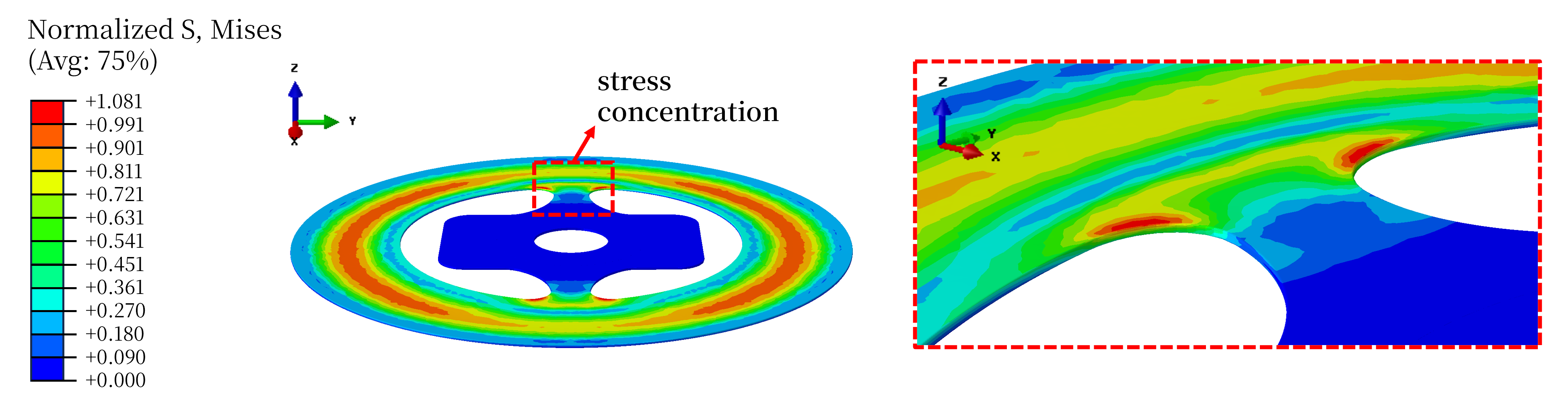}
        \caption{Normalized von Mises stress contour and deformed configuration under the bulging load case.}
        \label{fig:bulging_contour}
    \end{subfigure}
    \vspace{0.5mm}
    \caption{Normalized von Mises stress contours and deformed configurations of the initial nominal design under the bending and bulging load cases. For each load case, the stress values are normalized by the corresponding allowable CVaR limit used in the optimization. The highlighted regions indicate the primary stress-concentration locations.}
    \label{fig:stress_contour}
\end{figure}

The bending and bulging responses were evaluated through separate explicit analyses using the pressure loads shown in Figures~\ref{fig:bending} and~\ref{fig:bulging}, respectively. We defined general contact between opposing valve and piston surfaces that were initially in contact or could come into contact during deformation. The contact formulation used hard contact in the normal direction and frictionless behavior in the tangential direction. We prescribed the pressure magnitudes using the available operating and displacement data, including the known valve opening in the bending case, rather than deriving them from a fluid-dynamic analysis.
%
%
For each load case, we defined the stress response as the maximum von Mises stress across all valve integration points and throughout the analysis step. We retained stresses at contact edges and constrained regions and applied no nodal averaging. The averaging threshold shown in Figure~\ref{fig:stress_contour} affects only the contour visualization and not the extracted stress response. In the bending analysis, we computed the valve opening height as the maximum $z$-direction displacement of the valve nodes relative to the initial configuration.

\begin{figure}[htbp]
    \centering
    \includegraphics[width=0.6\linewidth]{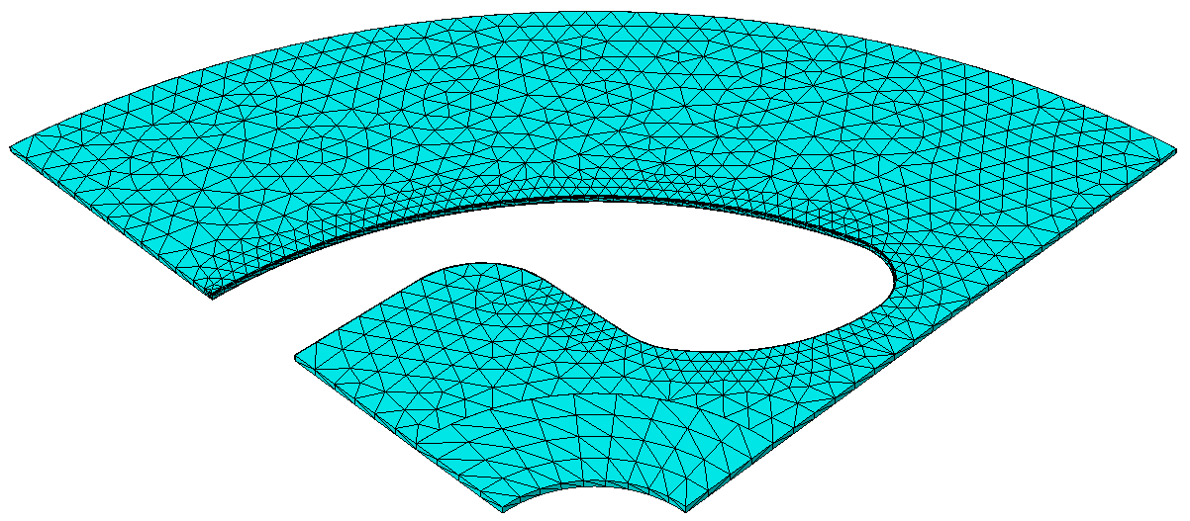}
    \caption{Tetrahedral finite-element discretization of the valve quarter model.}
    \label{fig:mesh}
\end{figure}

We analyzed the model in Abaqus/CAE 2025 using a Dynamic, Explicit step. The Freeflex valve had Young's modulus $210~\mathrm{GPa}$, Poisson's ratio $0.29$, and density $7.7~\mathrm{g/cm^{3}}$; the aluminum piston used $68.9~\mathrm{GPa}$, $0.33$, and $2.7~\mathrm{g/cm^{3}}$, respectively. We defined plasticity by the yield-stress--plastic-strain pairs $(1600~\mathrm{MPa},0)$ and $(2020~\mathrm{MPa},0.07)$ for the valve and $(276~\mathrm{MPa},0)$ and $(310~\mathrm{MPa},0.17)$ for the piston.

We discretized both components using ten-node modified quadratic tetrahedral elements (C3D10M). As shown in Figure~\ref{fig:mesh}, we used a global mesh size of $0.5$ and refined the stress-concentration regions of the valve using local seeds of $0.25$. The resulting model comprised $78{,}924$ elements and $119{,}838$ nodes. We selected this mesh resolution by progressively refining the mesh and monitoring the convergence of the stress responses while accounting for the associated computational cost.

We set the total step time to $0.22~\mathrm{s}$. For the bending analysis, we used a tabular amplitude to maintain the prescribed bending pressure $p_{\mathrm{bending}}$ from the beginning of the step to $0.10~\mathrm{s}$ and then reduce it to zero by $0.11~\mathrm{s}$. For the bulging analysis, we increased the pressure linearly from zero to the prescribed bulging pressure $p_{\mathrm{bulging}}$ during the first $0.10~\mathrm{s}$ and then maintained this value for the remainder of the step. We applied mass scaling to the entire model with a target stable time increment of $1.0\times10^{-6}~\mathrm{s}$. The procedure scaled the element mass only when the stable time increment fell below this target. To generate the dataset, we sampled the full design domain using Latin hypercube sampling and created the corresponding parameterized CAD geometries in CATIA V5. We automated the analysis workflow using the Abaqus Macro Manager script that sequentially replaced the CAD geometry, constructed and solved the finite-element model, and extracted the bending-stress, bulging-stress, and valve-opening responses. Each analysis used 32 CPU cores. Figure~\ref{fig:stress_contour} shows the normalized bending and bulging stress contours at the initial nominal design.

We constructed separate DD-GPCE surrogate models for the bending stress, bulging stress, and valve opening height. Each model used $N=24$ random variables, a maximum interaction order of $S=1$, and a polynomial order of $m=3$, resulting in $L_{24,1,3}=73$ basis functions. The initial global training datasets contained 828, 810, and 394 HF samples for the bending stress, bulging stress, and opening height, respectively. For uncertainty propagation, we evaluated each surrogate using $L=10{,}000$ random samples. A preliminary convergence assessment showed that the estimated CVaR values stabilized beyond this sample size. We evaluated the bending and bulging stress constraints at a risk level of $\beta=0.99$ and set the confidence-interval parameter and Tikhonov regularization parameters to $\alpha=0.5$ and $\lambda=10^{-6}$, respectively. At each tail-correction step, we selected $N_{1}=3$ samples during the predictive-variance-based exploration stage and $N_{2}=3$ samples during the subsequent prediction-shift-based exploitation stage.

We performed the optimization using the DE algorithm with the \texttt{best1bin} strategy. The mutation factor was varied within $(0.0,1.5)$, the recombination constant was set to $0.7$, and the population-size factor and convergence tolerance were specified as $5$ and $0.01$, respectively. We fixed the random seed at 42, disabled local polishing, and configured immediate population updates with a single worker. At the end of each generation, we corrected and updated the surrogate models. We then appended the HF samples acquired for tail-region correction to the global training datasets and reused them in subsequent generations. The optimization terminated after 25 generations.

\subsubsection{Optimization results}

Table~\ref{tab:design_variables} reports normalized changes from $\mathbf{d}^{(0)}$ to $\mathbf{d}^{\ast}$; unity denotes the initial value. Optimization reduces mass by 1.05\%. The largest changes occur in $d_2$, $d_3$, $d_5$, and $d_{15}$--$d_{17}$, with $d_3$ near its lower bound and $d_{15}$ near its upper bound. Without a separate sensitivity analysis, these are coupled optimizer-selected changes, not independent importance measures.

Table~\ref{tab:optimization_results} uses independent post-optimization samples. At each design, more than 300 HF samples trained separate DD-GPCE models, and six additional HF samples corrected the bending and bulging tails. The values are therefore independent MF tail-corrected estimates, not crude-MCS references. Reported HF counts cover global training and corrections during optimization but exclude this post-optimization analysis.

We normalize each stress by its allowable CVaR, so unity is the constraint threshold. The initial bending and bulging CVaR values, 1.1791 and 1.2649, both exceed their limits. Figure~\ref{fig:stress_contour} shows that bending concentrates stress near the curved transition from the flexible valve to the supported region. Bulging raises stress across the contact region but reaches its maximum at the same transition; this region therefore governs both load cases.

\begin{figure}[htbp]
    \centering
    \begin{subfigure}[b]{0.35\textwidth}
        \centering
        \includegraphics[width=\textwidth]{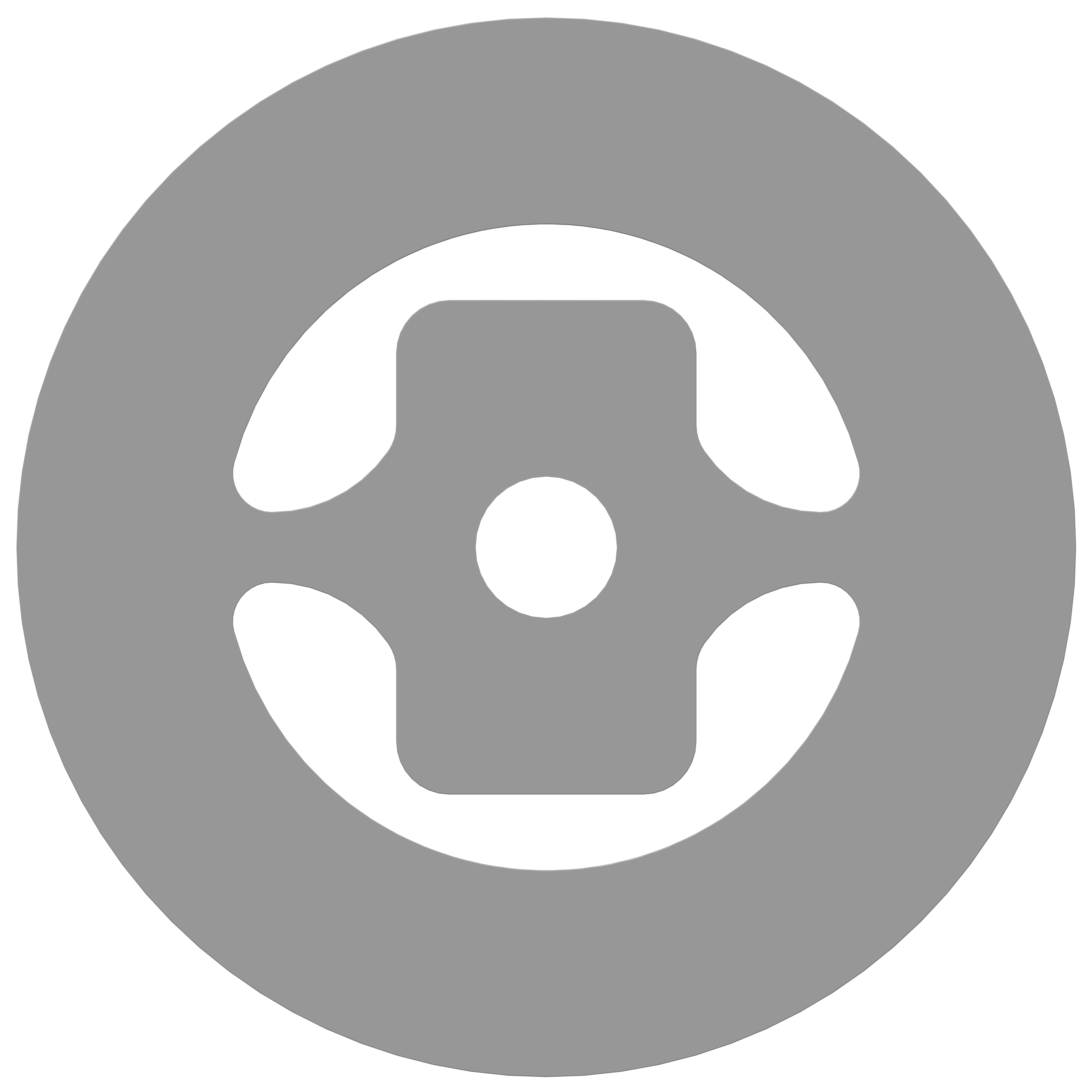}
        \caption{Initial design}
        \label{fig:init_design}
    \end{subfigure}
    \qquad
    \begin{subfigure}[b]{0.35\textwidth}
        \centering
        \includegraphics[width=\textwidth]{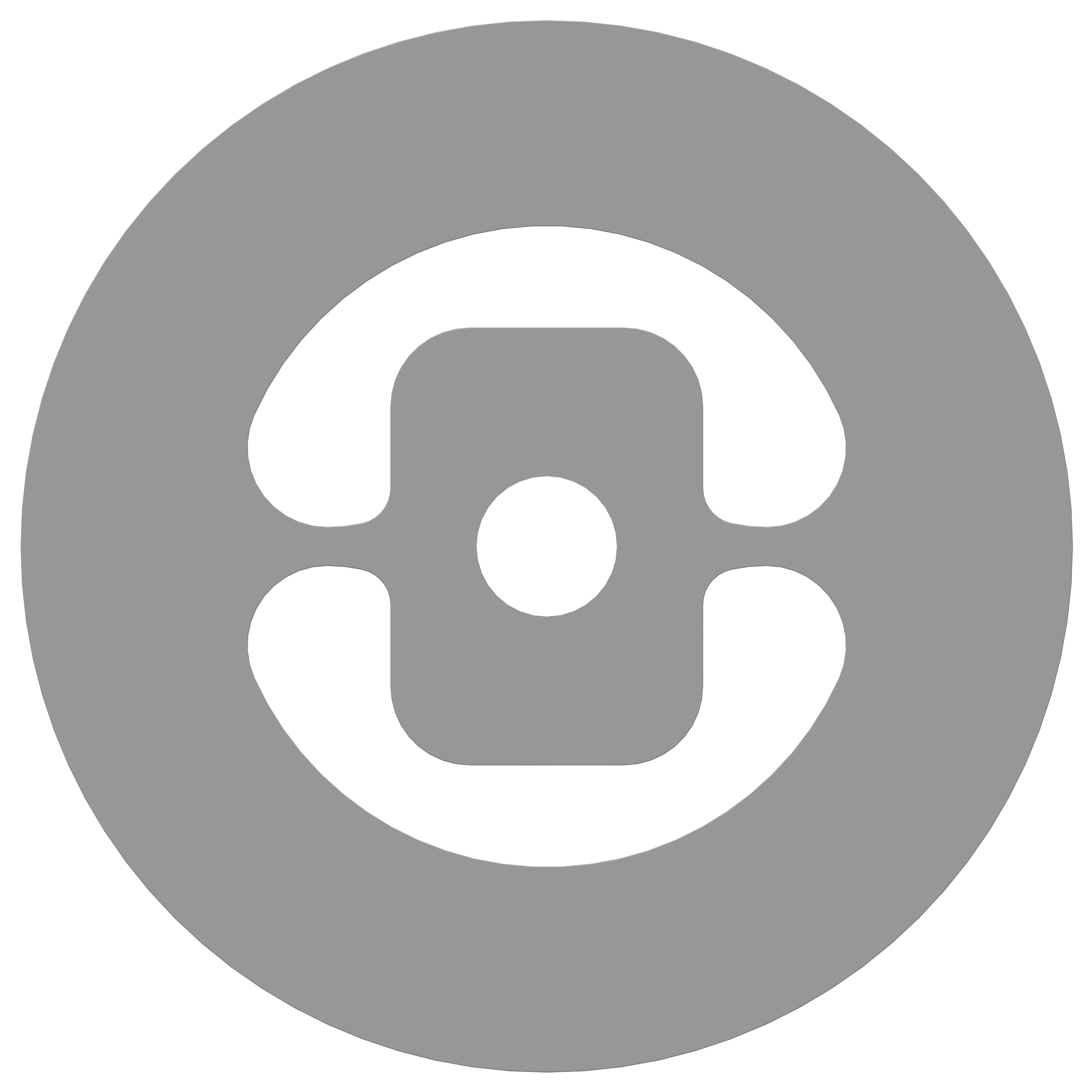}
        \caption{Optimal design}
        \label{fig:opt_design}
    \end{subfigure}
    \caption{Geometric comparison of the suction valve at the initial and optimal designs.}
    \label{fig:design_change}
\end{figure}

Figure~\ref{fig:design_change} compares the initial and optimal valve geometries and shows a pronounced modification of this stress-critical transition region. The optimal design smooths the transition by reducing the local curvature and increasing the radius of curvature. This modification transfers the bending deformation more gradually and thereby reduces the stress concentration under both loading conditions.

At the optimal design, the normalized bending and bulging CVaR values decrease to $0.9706$ and $0.9926$, respectively. Both values satisfy their risk constraints, with normalized feasibility margins of $0.0294$ and $0.0074$. The mean valve opening height, normalized by the prescribed minimum requirement $h_{\min}$, decreases from $1.0830$ to $1.0243$ but remains above the lower bound of unity. Among the three constraints, the bulging-stress constraint lies closest to its active limit.

\begin{figure}[!ht]
    \centering

    \begin{subfigure}[t]{0.7\linewidth}
        \centering
        \includegraphics[width=\linewidth]
        {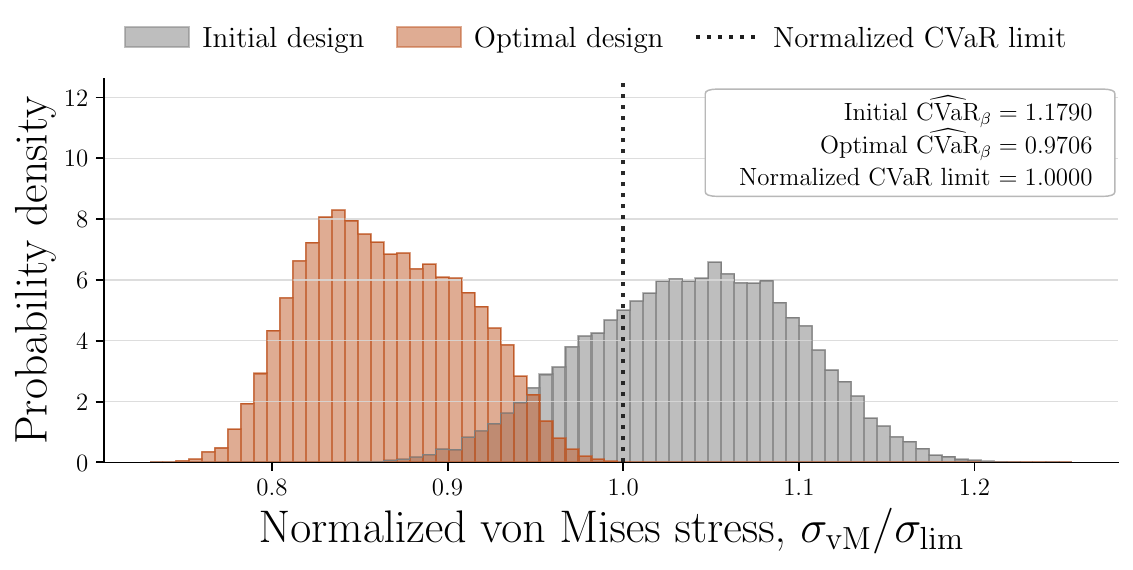}
        \caption{Normalized bending-stress distributions at the initial and optimal designs.}
        \label{fig:bending_dist}
    \end{subfigure}

    \vspace{5mm}

    \begin{subfigure}[t]{0.7\linewidth}
        \centering
        \includegraphics[width=\linewidth]
        {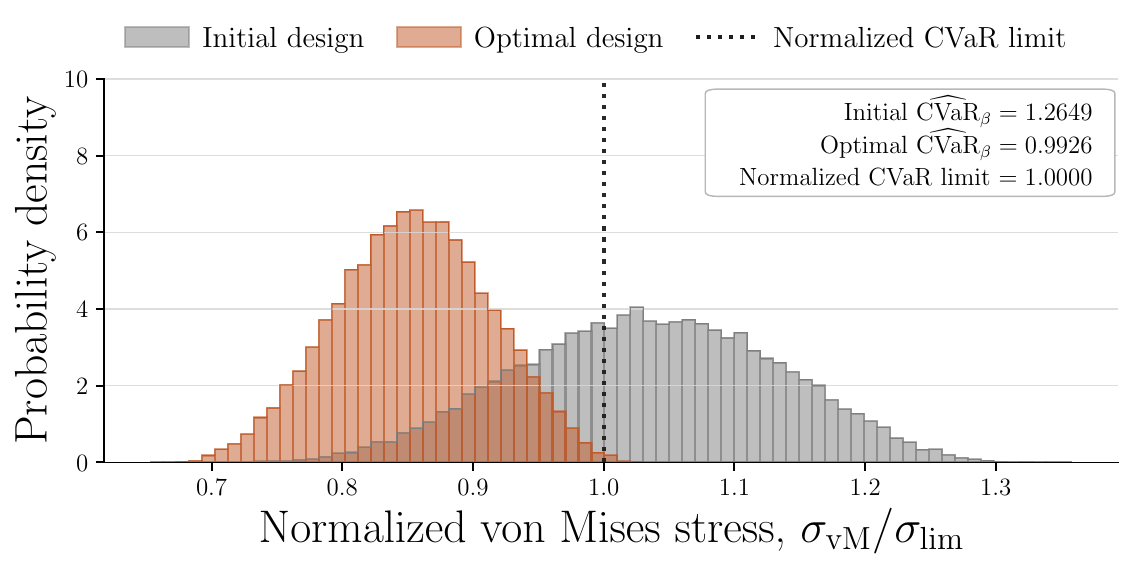}
        \caption{Normalized bulging-stress distributions at the initial and optimal designs.}
        \label{fig:bulging_dist}
    \end{subfigure}

    \caption{Normalized von Mises stress distributions estimated using the MF tail-correction models at the initial and optimal designs: (a) bending and (b) bulging. For each load case, the stress values are normalized by the corresponding allowable CVaR limit. The dotted vertical lines indicate the normalized CVaR limit of unity.}
    \label{fig:stress_distribution}
\end{figure}

Figure~\ref{fig:stress_distribution} shows that optimization shifts both normalized stress distributions leftward. We use the MF-tail-corrected models with $L=10{,}000$ samples and a common random seed at both designs. For bending, the mean, standard deviation, and CVaR decrease by 17.02\%, 29.49\%, and 17.68\%; for bulging, they decrease by 17.36\%, 40.12\%, and 21.53\%. Thus, lower tail risk reflects reductions in both location and dispersion. The vertical lines mark the normalized CVaR limits, not pointwise stress limits.

The modest mass reduction is consistent with the small final feasibility margins, especially for bulging stress. The workflow uses 1,934 finite-element jobs and $269.2~\mathrm{h}$ of CPU time, including global training and optimization. This example demonstrates the applicability of the proposed method to nonlinear contact under a limited HF budget.

    





\begin{table}[htb]
    \centering
    \begin{threeparttable}
        \caption{Optimization results for the suction-system design}
        \label{tab:optimization_results}

        \scriptsize
        \renewcommand{\arraystretch}{1.12}
        \setlength{\tabcolsep}{8pt}

        \begin{tabular}{@{}lllllcc@{}}
            \toprule
            Quantity of interest
            & Statistics
            & Initial
            & Optimal
            & Variation
            & No. of HF evaluations
            & CPU time \\
            \midrule

            Bending stress\tnote{a}
            & Mean [-]
            & 1.0425
            & 0.8651
            & $-17.02\%$
            & \multirow[c]{3}{*}{909}
            & \multirow[c]{9}{*}{$269.2~\mathrm{h}$} \\

            & Std. Dev. [-]
            & 0.0626
            & 0.0442
            & $-29.49\%$
            &
            & \\

            & $\mathrm{CVaR}_{\beta}$ [-]
            & 1.1791
            & 0.9706
            & $-17.68\%$
            &
            & \\

            \addlinespace[4pt]

            Bulging stress\tnote{a}
            & Mean [-]
            & 1.0320
            & 0.8528
            & $-17.36\%$
            & \multirow[c]{3}{*}{1025}
            & \\

            & Std. Dev. [-]
            & 0.0932
            & 0.0558
            & $-40.12\%$
            &
            & \\

            & $\mathrm{CVaR}_{\beta}$ [-]
            & 1.2649
            & 0.9926
            & $-21.53\%$
            &
            & \\

            \addlinespace[4pt]

            Valve opening height\tnote{b}
            & Mean [-]
            & 1.0830
            & 1.0243
            & $-5.42\%$
            & \multirow[c]{2}{*}{--}
            & \\

            & Std. Dev. [-]
            & 0.0575
            & 0.0559
            & $-2.87\%$
            &
            & \\

            \addlinespace[4pt]

            System mass\tnote{c}
            & Deterministic [-]
            & 1.0000
            & 0.9895
            & $-1.05\%$
            &
            & \\

            \bottomrule
        \end{tabular}

        \begin{tablenotes}[flushleft]
            \footnotesize
            \item[a] The bending and bulging stress statistics are normalized
            by their corresponding allowable CVaR limits,
            $\sigma_{\mathrm{lim}}^{(\mathrm{bending})}$ and
            $\sigma_{\mathrm{lim}}^{(\mathrm{bulging})}$, respectively.
            A normalized CVaR value of unity represents the corresponding
            constraint limit.
            \item[b] The valve-opening statistics are normalized by the
            prescribed minimum opening requirement $h_{\min}$. A normalized
            mean value of unity represents the lower constraint bound. The valve-opening responses were extracted from the same high-fidelity bending simulations; therefore, no separate high-fidelity evaluations were required for the valve-opening response.
            \item[c] The system mass is normalized by the initial system mass
            $m_0$. Accordingly, unity represents the mass of the initial
            design.
        \end{tablenotes}
    \end{threeparttable}
\end{table}
\FloatBarrier

\section{Conclusions} \label{sec:Conclusions}

This study developed an MF tail-correction method for risk-averse design optimization under CVaR constraints. The method used DD-GPCE as a global surrogate and directed limited HF evaluations to the response region that most strongly influenced the CVaR estimate. We quantified the finite-sample prediction uncertainty of DD-GPCE using the covariance of the estimated expansion coefficients and used this uncertainty to construct a confidence-interval-based $\epsilon$-risk region. Within this region, a two-stage adaptive strategy selected HF samples based on predictive variance and weighted prediction shift. We then represented the local surrogate discrepancy using a Tikhonov-regularized residual expansion to provide the tail-localized additive correction.


For the Griewank function, MF tail correction achieved lower MRD and N-RMSD values than MF importance sampling and standard DD-GPCE using four additional HF samples. It also converged near the theoretical optimum with fewer HF evaluations. For the ten-bar truss, the method produced a feasible design close to the crude-MCS reference using 618 HF evaluations, whereas the considered MF importance sampling configurations required more evaluations and retained small constraint violations.

The suction valve example demonstrated applicability to a nonlinear finite-element problem with 24 random variables and 18 design variables. The optimized design reduced the mass by 1.05\% and the bending and bulging CVaR values by 17.68\% and 21.53\%, respectively, while satisfying all constraints. The procedure required 1,934 finite-element analyses and $269.2~\mathrm{h}$ of CPU time. Because this example did not include crude-MCS or MF importance sampling comparisons, it demonstrates applicability but does not establish a direct computational-efficiency advantage.

The current confidence-interval model assumes homoscedastic Gaussian regression residuals and may not capture spatially varying prediction errors. Future work will develop heteroscedastic uncertainty models, incorporate epistemic uncertainty alongside aleatoric uncertainty in the optimization, and quantify confidence in the resulting optimal design.

\section*{Acknowledgments}
This work was partially supported by the National Research Foundation of Korea (NRF) grant funded by the Korea government (MSIT) (No. RS-2025-00560781) and partially by LG Electronics. The authors thank Jongeun Oh, Youngkyun Lim, and Jongtae Her of LG Electronics for their helpful discussions on the application of the suction valve system.

\bibliography{paper_template_RED/references.bib}
\bibliographystyle{unsrt}

\begin{appendix}

\section{Generalized polynomial chaos expansion}
\label{app:gpce}

The generalized polynomial chaos expansion (GPCE) represents a square-integrable random output in terms of multivariate orthonormal polynomials of the input random variables. When the components of $\mathbf{X}=(X_{1},\ldots,X_{N})^{\top}$ are statistically dependent, their joint probability measure is generally non-product-type. Therefore, measure-consistent multivariate orthonormal polynomials cannot, in general, be constructed as tensor products of univariate orthonormal polynomials. Instead, they are constructed directly from the joint probability measure $f_{\mathbf{X}}(\mathbf{x})\,d\mathbf{x}$.

Let $\mathbf{j}=(j_{1},\ldots,j_{N})\in\mathbb{N}_{0}^{N}$ be an $N$-dimensional multi-index. For $\mathbf{x}=(x_{1},\ldots,x_{N})^{\top}\in\mathcal{A}^{N}$, the corresponding monomial is
$\mathbf{x}^{\mathbf{j}}=x_{1}^{j_{1}}\cdots x_{N}^{j_{N}}$, with total degree
$|\mathbf{j}|=j_{1}+\cdots+j_{N}$. For a prescribed polynomial order $m\in\mathbb{N}_{0}$, define the total-degree multi-index set
\begin{equation}
    \mathcal{J}_{m}
    :=
    \left\{
        \mathbf{j}\in\mathbb{N}_{0}^{N}
        :
        |\mathbf{j}|\leq m
    \right\},
    \label{eq:app_gpce_index}
\end{equation}
whose cardinality is
\begin{equation}
    L_{N,m}
    :=
    |\mathcal{J}_{m}|
    =
    \sum_{q=0}^{m}
    \binom{N+q-1}{q}
    =
    \binom{N+m}{m}.
    \label{eq:app_gpce_cardinality}
\end{equation}

Let
\[
    \boldsymbol{\Psi}_{m}(\mathbf{x})
    =
    \left(
        \Psi_{1}(\mathbf{x}),
        \ldots,
        \Psi_{L_{N,m}}(\mathbf{x})
    \right)^{\top}
\]
denote the multivariate orthonormal polynomial vector consistent with
$f_{\mathbf{X}}(\mathbf{x})\,d\mathbf{x}$. Then, an output random variable
$y(\mathbf{X})\in L^{2}(\Omega,\mathcal{F},\mathbb{P})$ can be approximated by the $m$th-order GPCE
\begin{equation}
    \widetilde{y}_{m}(\mathbf{X})
    =
    \sum_{i=1}^{L_{N,m}}
        c_{i}\Psi_{i}(\mathbf{X}),
    \label{eq:app_regular_gpce}
\end{equation}
where the expansion coefficients are defined by
\begin{equation}
    c_{i}
    :=
    \mathbb{E}
    \left[
        y(\mathbf{X})\Psi_{i}(\mathbf{X})
    \right]
    =
    \int_{\mathcal{A}^{N}}
        y(\mathbf{x})
        \Psi_{i}(\mathbf{x})
        f_{\mathbf{X}}(\mathbf{x})
        \,d\mathbf{x},
    \qquad
    i=1,\ldots,L_{N,m}.
    \label{eq:app_gpce_coefficient}
\end{equation}

The GPCE retains all polynomial terms with total degree not exceeding $m$, whereas DD-GPCE reduces the basis by limiting the maximum interaction order among the input variables. The construction of the corresponding measure-consistent orthonormal polynomial basis is described in~\ref{app:build_basis}. We refer to~\cite{lee2020practical} for more details.

\section{Three steps for measure-consistent orthonormal polynomials for DD-GPCE}
\label{app:build_basis}

Let $\mathcal{I}_{N}:=\{1,\ldots,N\}$ be the complete input-index set. For an $N$-dimensional multi-index
$\mathbf{j}=(j_{1},\ldots,j_{N})\in\mathbb{N}_{0}^{N}$, let
$|\mathbf{j}|=\sum_{k=1}^{N}j_k$ denote its total degree and
$\|\mathbf{j}\|_{0}$ the number of its nonzero components. For a prescribed maximum interaction order $S$ and polynomial order $m$, the reduced multi-index set used in DD-GPCE is defined as
\begin{equation}
    \mathcal{J}_{S,m}
    :=
    \left\{
        \mathbf{j}\in\mathbb{N}_{0}^{N}
        :
        |\mathbf{j}|\leq m,\;
        \|\mathbf{j}\|_{0}\leq S
    \right\}.
    \label{eq:app_reduced_index}
\end{equation}
This set retains polynomial terms involving at most $S$ input variables and has cardinality
\begin{equation}
    L_{N,S,m}
    =
    1+
    \sum_{s=1}^{S}
    \binom{N}{s}
    \binom{m}{s}.
    \label{eq:app_reduced_cardinality}
\end{equation}
Let the elements of $\mathcal{J}_{S,m}$ be arranged as
$\mathbf{j}^{(1)},\ldots,\mathbf{j}^{(L_{N,S,m})}$, with
$\mathbf{j}^{(1)}=\mathbf{0}$. The measure-consistent orthonormal polynomial basis is then constructed through the following three steps.

\begin{enumerate}[label=\textbf{Step \arabic*}:, leftmargin=*]
    \item \textbf{Construct the reduced monomial vector.}

    For $\mathbf{x}=(x_1,\ldots,x_N)^{\top}$, define
    \begin{equation}
        \mathbf{M}_{S,m}(\mathbf{x})
        :=
        \left(
            \mathbf{x}^{\mathbf{j}^{(1)}},
            \ldots,
            \mathbf{x}^{\mathbf{j}^{(L_{N,S,m})}}
        \right)^{\top},
        \qquad
        \mathbf{x}^{\mathbf{j}}
        =
        \prod_{k=1}^{N}x_k^{j_k}.
        \label{eq:app_monomial_vector}
    \end{equation}

    \item \textbf{Construct the monomial moment matrix.}

    The moment matrix associated with the joint probability measure of $\mathbf{X}$ is
    \begin{equation}
        \mathbf{G}_{S,m}
        :=
        \mathbb{E}
        \left[
            \mathbf{M}_{S,m}(\mathbf{X})
            \mathbf{M}_{S,m}^{\top}(\mathbf{X})
        \right]
        =
        \int_{\mathcal{A}^{N}}
            \mathbf{M}_{S,m}(\mathbf{x})
            \mathbf{M}_{S,m}^{\top}(\mathbf{x})
            f_{\mathbf{X}}(\mathbf{x})
            \,d\mathbf{x}.
        \label{eq:app_moment_matrix}
    \end{equation}
    For a general joint probability density, this matrix can be estimated using numerical integration or sampling~\cite{lee2020practical}.

    \item \textbf{Apply the whitening transformation.}

    Let $\mathbf{W}_{S,m}$ be a whitening matrix obtained from the Cholesky factorization of the symmetric positive-definite moment matrix~\cite{rahman2018polynomial}, such that
    \begin{equation}
        \mathbf{G}_{S,m}
        =
        \mathbf{W}_{S,m}^{-1}
        \mathbf{W}_{S,m}^{-\top}.
        \label{eq:app_whitening_relation}
    \end{equation}
    The DD-GPCE basis vector is then generated as
    \begin{equation}
        \boldsymbol{\Psi}_{S,m}(\mathbf{x})
        =
        \mathbf{W}_{S,m}
        \mathbf{M}_{S,m}(\mathbf{x}).
        \label{eq:app_basis_whitening}
    \end{equation}
    Accordingly,
    \begin{equation}
        \mathbb{E}
        \left[
            \boldsymbol{\Psi}_{S,m}(\mathbf{X})
            \boldsymbol{\Psi}_{S,m}^{\top}(\mathbf{X})
        \right]
        =
        \mathbf{I}_{L_{N,S,m}},
        \label{eq:app_basis_orthonormality}
    \end{equation}
    so that the resulting multivariate polynomials are orthonormal with respect to the complete joint probability measure of $\mathbf{X}$. Statistical dependence among the inputs is therefore incorporated directly through the monomial moment matrix.

\end{enumerate}
The resulting $S$-variate, $m$th-order DD-GPCE approximation can therefore be expressed as
\begin{equation}
    \widetilde{y}_{S,m}(\mathbf{X})
    =
    \mathbf{c}^{\top}
    \boldsymbol{\Psi}_{S,m}(\mathbf{X})
    =
    \sum_{i=1}^{L_{N,S,m}}
        c_{i}\Psi_{i}(\mathbf{X}).
    \label{eq:app_ddgpce_final_expansion}
\end{equation}
Although the reduced monomial set omits interactions involving more than $S$ variables, the resulting basis remains orthonormal with respect to the complete joint probability measure of $\mathbf{X}$. Therefore, statistical dependence among the input variables is incorporated through the monomial moment matrix rather than through a tensor product of marginal polynomial bases. When the moment matrix is approximated numerically, the resulting basis satisfies the orthonormality condition up to the numerical integration or sampling error.

\end{appendix}
\end{document}